\documentclass[11pt,reqno]{amsart}

\usepackage{amsmath}
\usepackage{amsxtra}
\usepackage{amscd}
\usepackage{amsthm}
\usepackage{amsfonts}
\usepackage{amssymb}
\usepackage{eucal}
\usepackage[all]{xy}

\numberwithin{equation}{section}

\newtheorem{theorem}{Theorem}[subsection]
\newtheorem{proposition}[theorem]{Proposition}
\newtheorem{lemma}[theorem]{Lemma}
\newtheorem{corollary}[theorem]{Corollary}
\theoremstyle{definition}

\newtheorem{remark}[theorem]{Remark}
\newtheorem{example}[theorem]{Example}

\newcommand{\nc}{\newcommand}
\nc{\C}{{\mathbb C}}
\def\m{{\bf m}}
\let\al\alpha
\newcommand{\half}{\frac12}
\newcommand{\Z}{{\mathbb Z}}

\newcommand{\Ref}[1]{{$($\ref{#1}$)$}}
\newcommand{\bean}{\begin{eqnarray*}}
\newcommand{\eean}{\end{eqnarray*}}
\newcommand{\be}{$$}
\newcommand{\bea}{\begin{eqnarray}}
\newcommand{\ena}{\end{eqnarray}}
\newcommand{\ee}{$$}
\def\ov[#1,#2]{\overset{\scriptstyle #1}{#2}}
\nc{\verm}{M_{k,l}}
\def\wsl{\widehat{\mathfrak{sl}}_2}
\nc{\slt}{{\widehat{\mathfrak{sl}}_2}}
\nc{\inte}{{L_{k,l}}}
\nc{\pkln}{{\mathcal P}_{k,l}^{(N)}}
\nc{\ckln}{{c_{k,l}^{(N)}}}
\def\qbin[#1;#2]{{\left[\matrix{\displaystyle #1}\\{\displaystyle #2}\endmatrix\right]}}
\def\geq{\ge}
\def\leq{\le}

\def\cH{\widehat{\mathcal H}}

\nc{\Vk}{{\mathfrak V}_k}

\newcommand\no{\nonumber}

\def\xx[#1,#2]{{#1}^{(#2)}}
\def\yy[#1,#2,#3]{{#1}^{(#2)}_{#3}}

\newcommand{\res}{{\operatorname{res}}}

\def\v(#1;#2){\Bigl({#1\atop#2}\Bigr)}

\nc{\alb}{{\boldsymbol{\alpha}}}
\nc{\beb}{{\boldsymbol{\beta}}}
\def\ra{\rangle}
\def\la{\langle}
\def\hookdownarrow%
{\Big\downarrow\kern-.14cm\smash{\raise.4cm\hbox{$\scriptstyle\cap$}}}

\begin{document}
\title[Combinatorics of Coinvariants]{Combinatorics of the $\wsl$
Spaces of Coinvariants II}
\author{B. Feigin, R. Kedem,
 S. Loktev, T. Miwa
 and E. Mukhin}
\address{BF: Landau institute for Theoretical Physics, Chernogolovka,
142432, Russia}
\address{RK: Dept. of Mathematics, U. Massachusetts, Amherst, MA 01003.}
\address{SL: Independent University of Moscow, B. Vlasievsky per, 11,
Moscow}
\address{TM: RIMS, Kyoto University, Kyoto 606 Japan.}
\address{EM: Dept. of Mathematics, U. California, Berkeley, CA 94720.}

\begin{abstract}
The spaces of coinvariants are
quotient spaces of integrable $\widehat{\mathfrak{sl}}_2$
modules by subspaces generated by actions of certain subalgebras
labeled by a set of points on a complex line. When
all the points are distinct, the spaces of coinvariants
essentially coincide with the spaces of conformal blocks
in the WZW conformal field theory and their dimensions are 
given by the Verlinde rule.

We describe monomial bases for the $\wsl$ spaces of coinvariants, In
particular, we prove that the spaces of coinvariants have
the same dimensions when all the points coincide.
We establish recursive relations satisfied by the monomial bases and
the corresponding characters of the spaces of coinvariants.
For the proof we use filtrations of the $\widehat{\mathfrak{sl}}_2$ modules,
and further filtrations on the adjoint graded spaces for the first filtrations.

This paper is the continuation of \cite{FKLMM}.
\end{abstract}
\maketitle

\section{Introduction}\label{INTRO}
Let $e,h,f$ be the standard basis of the Lie algebra $\mathfrak{sl}_2$,
and $e(x)=\sum_{i\in\Z}e_ix^i$, $h(x)=\sum_{i\in\Z}h_ix^i$,
$f(x)=\sum_{i\in\Z}f_ix^i$ the currents of
$\wsl={\mathfrak sl}_2\otimes\C[t,t^{-1}]\oplus\C c\oplus\C d$.
By $\mathfrak{a}^{(M,N)}$ we denote the Lie subalgebra of
$\wsl$ generated by $e_i$ ($i\geq M$) and $f_i$ ($i\geq N$),
where $M,N$ are non-negative integers.
Let $L_{k,l}$ be the integrable irreducible representation of
$\wsl$ of level $k$ and highest weight $l\in\{0,\dots,k\}$.
We suppose that the highest weight vector of $L_{k,l}$ is annihilated by
$e_i$ ($i\leq0$), $h_j$ and $f_j$ ($j<0$).

In this paper, we study the space of coinvariants
$$
L^{(M,N)}_{k,l}=L_{k,l}/\mathfrak{a}^{(M,N)}L_{k,l}.
$$
We prove that $L^{(M,N)}_{k,l}$ is a finite-dimensional space and its
dimension  is given by the Verlinde formula. Namely, we have
the recurrence formula
\begin{equation}
{\rm dim}\,L^{(M,N)}_{k,l}=\sum_{l',l''}{\rm dim}\,L^{(M,N-1)}_{k,l'},
\label{2}
\end{equation}
where the summation is over the pairs of non-negative
integers $(l',l'')$ such that $l$ has the same parity as $l'+l''$ and
$$
|l'-l''|\leq l\leq{\rm min}(l'+l'',2k-l'-l'').
$$

We prove (\ref{2}) in the following way.

Recall the Verlinde algebra for the $\widehat{\mathfrak sl}_2$-theory
of level $k$. It is a commutative algebra with the basis
$\{\pi_0,\dots,\pi_k\}$. The defining relations are
$$
\pi_{l''}\cdot\pi_{l'}=
\sum_{l\equiv l'+l''\bmod 2\atop|l'-l''|\leq l\leq{\rm min}(l'+l'',2k-l'-l'')}
\pi_l.
$$

First, we prove that
${\rm dim}\,L^{(M,N)}_{k,l}$ is greater than or equal to the number
$d^{(M+N)}_{k,l}$ defined by the formula
\begin{equation}
(\pi_0+\dots+\pi_k)^N=\sum_{0\leq l\leq k}d^{(N)}_{k,l}\pi_l.
\label{3}
\end{equation}
More precisely, we have the following result.
Consider a pair of polynomials $P(t), Q(t)$ with degrees ${\rm deg}\,P=M$
and ${\rm deg}\,Q=N$. Let $\mathfrak{a}(P,Q)$ be the subalgebra of
$\wsl$ generated by the subspaces $\C e\otimes P(t)\C[t]$
and $\C f\otimes Q(t)\C[t]$. It is clear that $\mathfrak{a}(t^M,t^N)$
is exactly $\mathfrak{a}^{(M,N)}$. Consider the family of
spaces of coinvariants:
$$
L_{k,l}(P,Q)=L_{k,l}/\mathfrak{a}(P,Q)L_{k,l}.
$$
It is well-known (for the proof see, e.g., the appendix of \cite{FKLMM})
that (\ref{2}) is true if we replace 
${\rm dim}\,L^{(M,N)}_{k,l}$ by ${\rm dim}\,L_{k,l}(P,Q)$ for ``generic''
polynomials $P,Q$. The word ``generic'' means that all $M+N$ zeros
of the polynomials $P,Q$ are distinct. It is easy to see that for arbitrary
$P,Q$ the dimension of $L_{k,l}(P,Q)$ is greater than or equal to the
dimension at a ``generic'' point. Formula (\ref{2}) actually means that the
family of vector spaces $L_{k,l}(P,Q)$ forms a vector bundle over the space
of parameters $(P,Q)$. In other words, the dimension does not ``jump''
at special points.

Let us try to understand why the formula (\ref{2})
is valid. The reason may be that there exists some natural filtration
in $L^{(M,N)}_{k,l}$ such that the corresponding graded space is isomorphic to
the direct sum  $\oplus_{l',l''}L^{(M,N-1)}_{k,l'}$. We do not know such
a filtration, nevertheless what we show in this paper is rather similar.

We replace $\wsl$ by the Lie algebra
$\cH={\mathcal H}\otimes\C[t,t^{-1}]$, where ${\mathcal H}$ is the
three-dimensional Heisenberg algebra with the basis $\{\bar e,\bar h,\bar f\}$
and the relations $[\bar e,\bar f]=\bar h$ and
$[\bar h,\bar e]=[\bar h,\bar f]=0$. As explained in \cite{FKLMM},
$\cH$ is the adjoint graded Lie algebra associated to $\wsl$
with respect to the filtration. Let
$F_0=\C c\oplus\C d$, $F_1=F_0\oplus\Bigl(\oplus_{i\in\Z}\C e_i\Bigr)
\oplus\Bigl(\oplus_{j\in\Z}\C f_j\Bigr)$ and $F_2=F_3=\dots=\wsl$.
The graded Lie algebra $F_0\oplus F_1/F_0\oplus F_2/F_1$ is isomorphic
to the sum $\cH\oplus\C c\oplus\C d$.

The filtration $\{F_i\}$ in $\wsl$ naturally induces
a filtration in any cyclic representation of $\wsl$.
Let $M$ be a cyclic representation and $v$ a cyclic vector in $M$.
We define the filtration $\{F_i(M)\}_{i\geq0}$ by the following inductive
procedure: $F_0(M)=\C v$, $F_{i+1}(M)=F_1\cdot F_i(M)$
($F_1\subset\wsl$). It is evident that the adjoint graded space
${\rm Gr}(M)=\oplus_iF_{i+1}(M)/F_i(M)$ is a cyclic representation of $\cH$.

We describe explicitly the $\cH$-module ${\rm Gr}(L_{k,l})$.

In the algebra $\cH$ we have the basis $\bar e_i,\bar h_i,\bar f_i$
($i\in\Z$) and the corresponding currents $\bar e(x),\bar h(x),\bar f(x)$.
Recall that in $L_{k,l}$ the following relations hold:
$e(x)^{k+1}=f(x)^{k+1}=0$. From this we deduce that in ${\rm Gr}(L_{k,l})$ 
$\bar e(x)^{k+1}=\bar f(x)^{k+1}=0$.
Define now the $\cH$-module $\bar L_{k,l}$ as follows. Let $\bar M$ be
the quotient of the module induced from the one-dimensional trivial
representation $\C\bar v$ of the subalgebra spanned by
$\bar e_i$ ($i\leq0$), $\bar h_j$ and $\bar f_j$ ($j<0$). Then, 
$\bar L_{k,l}$ is the quotient of  $\bar M$ 
by the submodule generated by the vectors $\bar e(x)^{k+1}w$,
$\bar f(x)^{k+1}w$ ($w\in \bar M$) and $\bar f_0^{l+1}\bar v$. We
prove the isomorphism of $\cH$-modules ${\rm Gr}(L_{k,l})\simeq\bar L_{k,l}$.

For simplicity let us assume that $M>0$ in this introduction.
In the main text, the case $M=0$ is also handled with a proper modification.
For the $\cH$-module $\bar L_{k,l}$ we define the spaces of coinvariants
$\bar L^{(M,N)}_{k,l}$ and show that 
${\rm dim}\,\bar L^{(M,N)}_{k,l}={\rm dim}\,L^{(M,N)}_{m,n}$.
Actually, we prove that (\ref{2}) is true for
${\rm dim}\,\bar L^{(M,N)}_{k,l}$. For this purpose we define a bigger set of
representations of $\cH$ depending on three parameters $l_1,l_2,l_3$.
For this set of $\cH$-modules, we define coinvariants and prove
an analogue of (\ref{2}) using suitable filtrations. Such recursive
description of the spaces $\bar L^{(M,N)}_{k,l}$ gives us
a collection of monomial bases in $\bar L^{(M,N)}_{k,l}$
(and therefore in $L^{(M,N)}_{k,l}$).

The plan of this paper is as follows. In Section 2, we give some
preliminary definitions and propositions. In Section 3, the bijection between
the combinatorial and Verlinde paths is given. In Section 4 several mappings
between Heisenberg modules are given. In Section 5, the recursion
relation for Heisenberg modules is proved. In Section 6, some miscellaneous
results are given. In Appendix, we give several lemmas for zero vectors
in Heisenberg modules.

\subsection{Acknowledgements} The research of B. Feigin was supported
by grants RFBR 99-01-01169 and INTAS-OPEN-97-1312.

\section{Preliminaries}\label{defs}
\subsection{Spaces of coinvariants}

Let $\mathfrak{sl}_2$ be the Lie algebra with basis $\{e,f,h\}$ such that
$$
[e,f]=h,\qquad [e,h]=-2e, \qquad [f,h]=2f.
$$
Let $(\;,\;)$ be the Killing form normalized by $(h,h)=2$.
The algebra $\slt=\mathfrak{sl}_2\otimes
\C[t,t^{-1}] \oplus \C c \oplus \C d$ is spanned by
$\{g_i=g\otimes t^i,c,d \ ;\ g=e,f,h
\in \mathfrak{sl}_2, i\in \Z\}$, where $c$ is a
central element and $d$ is the scaling operator, $[d,g_i]=i g_i$.
The commutators are given by (notice our normalization of the central
element) 
\bean 
[X\otimes f(t),Y\otimes g(t)]&=&[X,Y]\otimes
f(t)g(t)+c\cdot(X,Y)\res_{t=\infty}\frac{df}{dt}(t) g(t)\\
&=&[X,Y]\otimes
f(t) g(t)-c\cdot(X,Y)\res_{t=0}\frac{df}{d t}(t)
g(t).
\eean

Consider the Verma module $\verm=U(\slt) v_k[l]$, with $k\in\Z_{>0}$
and $l\in\{0,\dots,k\}$, generated by the highest weight vector
$v=v_k[l]$ satisfying
\bean
&&cv=kv,\quad h_0v=lv,\quad dv=0,\\
&&e_iv=0\ (i\leq 0),\quad h_iv=f_iv=0\ (i<0).
\eean
The integrable module
$\inte$ is the quotient of $\verm$ by the submodule generated by the vectors
$$
e(z)^{k+1}w\quad(\hbox{or equivalently $f(z)^{k+1}w$})
\quad\text{where $w\in\verm$}.
$$
Here we used the generating series notations $e(z)=\sum_{i\in\Z}
e_i z^i$, etc.. The module $\inte$ is an irreducible $\slt$ module.

\begin{remark}
Usually, the module $L_{k,l}$ is defined as the quotient of $M_{k,l}$ by
the relations $f_0^{l+1}v=e_1^{k-l+1}v=0$. These relations follow from the
relation $f(z)^{k+1}=0$. Indeed, 
the condition $f(z)^{k+1}v=0$ and $f_iv=0$ ($i<0$) implies $f_0^{k+1}v=0$.
Then, applying ${\rm ad}e_0$ to $f_0^{k+1}v=0$ and using the identity
$[e_0,f_0^{j+1}]v=(j+1)(l-j)f_0^jv$ and $h_0v=lv$, we obtain 
$f_0^{l+1}v=0$. Similarly, we have $e_1^{k-l+1}v=0$.
\end{remark}

The module $\inte$ is graded by $h_0$ and $d$. We denote its character
by $\chi_{k,l}$,
$$
\chi_{k,l}=\hbox{\rm ch}_{q,z}\inte
=\hbox{\rm trace}_{\inte}\left(q^dz^{h_0}\right).
$$

Let $M,N\in\Z_{\geq 0}$,
${\mathbf z}=(z_1,\ldots,z_M)\in\C^M$ and 
${\mathbf z}'=(z'_1,\ldots,z'_N)\in\C^N$, 
where $z_i, z_i'$ are not necessarily distinct.
Consider the following Lie subalgebra of $\slt$:
$$
\mathfrak{sl}_2^{(M,N)}({\mathbf z};{\mathbf z}')=
\Bigl[\C e\otimes\C[t]P(t)\Bigr]
\oplus\Bigl[\C h\otimes\C[t]P(t)Q(t)\Bigr]
\oplus\Bigl[\C f\otimes\C[t]Q(t)\Bigr],
$$
where $P(t)=\prod_{i=1}^M(t-z_i)$ and $Q(t)=\prod_{i=1}^N(t-z'_i)$.

The corresponding space of coinvariants is defined by
$$
L^{(M,N)}_{k,l}({\mathbf z};{\bf z}')=L_{k,l}/\mathfrak{sl}_2^{(M,N)}({\bf z};{\bf z}')L_{k,l}.
$$

The space $L^{(M,N)}_{k,l}({\bf z};{\bf z}')$ is finite-dimensional,
and for generic ${\bf z},{\bf z}'$ the
dimension is given in terms of the Verlinde algebra. Let us recall
its definition.

The level $k$ Verlinde algebra associated to $\mathfrak{sl}_2$
is a complex $k+1$-dimensional commutative associative algebra with
basis $\{\pi_l;\ 0\leq l \leq k\}$ and multiplication 
\bea\label{verlinde multiplication}
\pi_l\cdot\pi_{l'}=\sum_{\overset{\scriptstyle i=|l-l'|}{i+l-l':\ 
\hbox{\rm even}}}
^{\hbox{\rm min}(2k-l-l',l+l')}\pi_i.
\ena
For $N\in\Z_{>0}$, define the positive integer $d^{(N)}_{k,l}$ by
\bea\label{verlinde number}
(\pi_0+\pi_1+\cdots+\pi_k)^N=\sum_ld^{(N)}_{k,l}\pi_l.
\ena

\begin{theorem}\label{Ver}
$($\cite{FKLMM}$)$
If the points $z_1,\ldots,z_M,z_1',\dots,z_N'$ are distinct, then
$$
\hbox{\rm dim}\,L^{(M,N)}_{k,l}({\bf z};{\bf z}')=d^{(M+N)}_{k,l}.
$$
\end{theorem}

Theorem \ref{Ver} is proved in Appendix of \cite{FKLMM}.

We denote
$$
L^{(M,N)}_{k,l}=L^{(M,N)}_{k,l}((0,\dots,0),(0,\dots,0)).
$$
The space  $L^{(M,N)}_{k,l}$ inherits from $L_{k,l}$ a grading by 
$h_0$ and $d$.  We define the character 
\bea\label{coinv character def}
\chi_{k,l}^{(M,N)}=\hbox{\rm ch}_{q,z}L^{(M,N)}_{k,l}=\hbox{\rm trace}_{L^{(M,N)}_{k,l}}\left(q^dz^{h_0}\right).
\ena
We denote by $\Bigl(L^{(M,N)}_{k,l}\Bigr)_{s,e}$
the weight subspace of $L^{(M,N)}_{k,l}$ where
$h_0=s$ and $d=e$.
\subsection{A deformation argument}

We will need the following deformation lemma.

\begin{lemma}\label{deformation}$($\cite{FKLMM}$)$
Let $V$ be a vector space with filtration
$$
0=F_{-1}\subset F_0\subset F_1\subset\cdots\subset V, \qquad V=\bigcup_j F_j.
$$
Let $T_i:V\to V$, $(i\in I)$, be a set of linear maps with degree $d_i\geq0$,
i.e. $T_i(F_j)\subset F_{j+d_i}$.

Let ${\rm Gr}^F(V)=\oplus_j\;{\rm Gr}_j^F(V)=\oplus_j\; F_j/F_{j-1}$ and let  $\overline{T}_i$ be the induced graded maps. Then, we have the inequality for the dimensions
of the spaces of coinvariants:
\bea\label{deformation inequality}
{\rm dim}\;\left( V/\sum_iT_iV \right) \leq
{\rm dim}\;\left({\rm Gr}^F(V)/\sum_i\overline{T}_i\,{\rm Gr}^F(V)\right).
\ena
\end{lemma}
\begin{proof} See \cite{FKLMM}.
\end{proof}

\begin{remark}\label{rem}
Note that if in a particular case,
\Ref{deformation inequality} can be shown to be an equality,
and the image of $\sqcup_j\{v_\alpha^{(j)}\in F_j;1\leq\alpha\leq
{\rm dim}\,{\rm Gr}^F_j(V)/\sum_i\bar T_i{\rm Gr}^F_{j-d_i}(V)\}$ forms
a basis in ${\rm Gr}^F(V)/\sum_i\overline{T}_i{\rm Gr}^F(V)$, then the 
image of the same set of vectors forms a basis in $V/\sum_iT_iV$.
\end{remark}

\begin{proposition}\label{upper bound}
We have an inequality
$$
{\rm dim}\,
L_{k,l}^{(M,N)}({\bf z};{\bf z}')
\leq{\rm dim}\,L^{(M,N)}_{k,l}.
$$
\end{proposition}

\begin{proof}
Proposition \ref{upper bound} follows from Lemma \ref{deformation} applied to the vector space $L_{k,l}$ with the filtration induced by the $d$-grading,
and the set of operators
$$
\{e\otimes P(t)t^i, f\otimes Q(t)t^{i}, \;i\geq 0\}.
$$
\end{proof}
 
One of the main results which we prove in this paper is that in fact
\begin{equation}
{\rm dim}\,L_{k,l}^{(M,N)}({\bf z};{\bf z}')
={\rm dim}\,L^{(M,N)}_{k,l}
\end{equation}
(see Theorem \ref{MAINTH}).

\section{Combinatorial and Verlinde paths}\label{monomial}
We establish the bijection between two sets of ``paths'',
the combinatorial paths and the Verlinde paths. The former will
be used to describe the monomial basis of
$L_{k,l}^{(0,N)}({\bf z})$. The latter arises in the well-known
result in the $\wsl$-invariant conformal field theory.
The bijection leads to a recursive formula for the combinatorial paths.
This will be a key result in our proof of Theorem \ref{MAINTH}
given in Section 5.

\subsection{Combinatorial paths}
Let ${\mathcal C}_{k,l}^{(N)}$ be the set of all
$({\bf a};{\bf b})=(a_{N-1},\dots,a_0;b_{N-1},\dots,b_1\}
\in\{0,\dots,k\}^{2N-1}$ such that 
\bea\label{a0 condition}
a_0 \leq l, 
\ena
\bea\label{triangle condition}
a_i + b_{i+1} + a_{i+1}\leq k,\qquad b_i + a_i + b_{i+1}\leq k\qquad(i\geq 0),
\ena
\bea\label{trapezoid condition}
\sum_{s=i}^j b_s \leq k+\sum_{s=i+1}^{j-2} a_s\qquad(-1\leq i<j\leq\infty),
\ena
where we set $b_\infty=k, b_{-1}=l, b_0=a_{-1}=a_\infty=0$ and $a_i=b_i=0$ for $N\leq i <\infty$. 

Note, in particular, that
$$
b_{N-1}=0
$$
because of \Ref{trapezoid condition} with $i=N-1,j=\infty$.
        
In the above definition we assumed that $N$ is an positive (i.e., $N>0$)
integer.
For $N=1$ we have $a_0=l$ and $\sharp({\mathcal C}^{(1)}_{k,l})=1$.
There is a canonical injection
${\mathcal C}^{(N)}_{k,l}\rightarrow{\mathcal C}^{(N+1)}_{k,l}$.
We define
\begin{equation}
{\mathcal C}^{(0)}_{k,l}
=
\begin{cases}
{\mathcal C}^{(1)}_{k,0}&\text{ if $l=0$};\\
\emptyset&\text{ otherwise}.
\end{cases}
\label{Z}
\end{equation}

We call elements $({\bf a};{\bf b})\in{\mathcal C}_{k,l}^{(N)}$
combinatorial paths of length $N$, level $k$ and weight $l$.

It is convenient to write combinatorial paths
$({\bf a};{\bf b})\in {\mathcal C}_{k,l}^{(N)}$ as follows:
$$
\left( \begin{array}{llllllllllllll}
b_\infty= k  & \cdots    & 0 &    &b_{N-1} & \cdots&     & b_2  &     &  b_1 &   & 0 &  & b_{-1}=l\\
  \dots& 0 &   &  a_{N-1}&     & \cdots &a_2      & & a_1 &         & a_0 & & 0 &
\end{array}\right)
$$

One of the main purposes of this section is to prove that
the cardinality of the set ${\mathcal C}_{k,l}^{(N)}$ is given
by the Verlinde number $d_{k,l}^{(N)}$ defined by \Ref{3}:
$$
\sharp({\mathcal C}_{k,l}^{(N)})=d_{k,l}^{(N)}.
$$
As shown in the appendix of \cite{FKLMM},
we know that the Verlinde number $d_{k,l}^{(M+N)}$ gives the dimension
of the space $L_{k,l}^{(M,N)}({\bf z};{\bf z}')$ when $({\bf z};{\bf z}')$
are distinct.
Therefore, we can think of ${\mathcal C}_{k,l}^{(N)}$ as a set
which parametrizes a basis of coinvariants. In fact,
in Section \ref{SEC-COIN} we will prove that the images of vectors 
\bea
\{f_{N-1}^{a_{N-1}} h_{N-1}^{b_{N-1}}\dots
f_1^{a_1}h_1^{b_1}f_0^{a_0}v_k[l]\; ; \;({\bf a};{\bf b})\in{\mathcal C}^{(N)}_{k,l}\}
\ena
form a basis in $L_{k,l}^{(0,N)}({\bf z})$ for all ${\bf z}\in \C^N$
(Corollary \ref{FHBASIS}).

One can check that for $k=1$ the sets ${\mathcal C}^{(N)}_{1,0}$ and
${\mathcal C}^{(N)}_{1,1}$ coincide with sets ${\mathcal C}^{0,N}$ and
${\mathcal C}^{1,N}$ in \cite{FKLMM}. For $k=1$, Theorem \ref{MAINTH}
and Corollary \ref{FHBASIS} are proved in \cite{FKLMM}
(see Corollary 5.2.2 and Theorem 5.3.1 therein).

\bigskip 

\subsection{Verlinde paths}
In this section we recall the definition and some elementary properties of
Verlinde paths.

A triple $(\al,\beta,\gamma)\in\{0,\dots,k\}^3$ is called an admissible
Verlinde triple (of level $k$) if it can be written as
\begin{equation}
(\al,\beta,\gamma)\mapsto
\begin{pmatrix}\beta\\ \al\ \gamma\end{pmatrix}=
\begin{pmatrix} y+z \\ {x+y \ x+z}\end{pmatrix}
,\qquad x+y+z\leq k,\qquad x,y,z\in\Z_{\geq 0}.
\label{ADM}
\end{equation}
Then the numbers $x,y,z$ are determined uniquely,
\bea\label{xyz}
x(\al,\beta,\gamma)=\frac{\al+\gamma-\beta}{ 2},\qquad 
y(\al,\beta,\gamma)=\frac{\al+\beta-\gamma}{ 2},\qquad 
z(\al,\beta,\gamma)=\frac{\beta+\gamma-\al}{ 2}.
\ena

If $(\al,\beta,\gamma)$ is an admissible triple then any permuted triple,
e.g., $(\gamma,\al,\beta)$, is also admissible. Note, however, that formula
\Ref{xyz} is not symmetric in $\al,\beta,\gamma$.

Now, formula \Ref{verlinde multiplication}  for multiplication in the Verlinde algebra can be rewritten as
\bea\label{verlinde multiplication II}
\pi_l\cdot\pi_{l''}=\sum_{l':\,(l,l'',l')\;\; {\rm is}\;\;{\rm admissible}}\pi_{l'}.
\ena

A sequence
$(\alb;\beb)=(\al_1,\dots,\al_N;\beta_1,\dots\beta_{N-1})\in\{0,\dots,k\}^{2N-1}$
is called a Verlinde path of length $N$, weight $l$ and level $k$ if
$\al_1=l$ and $(\al_i,\beta_i,\al_{i+1})$, $i=1,\dots,N-1$ are
admissible Verlinde triples of level $k$. We set $\beta_N=\al_N$,
$\al_i=\beta_i=0$ for $i>N$, so that $(\al_i,\beta_i,\al_{i+1})$ form
Verlinde triples for all $i\in\Z_{>0}$.

We denote ${\mathcal P}_{k,l}^N$ the set of all Verlinde paths of length $N$, level $k$ and weight $l$. 

We define the following two maps associated with Verlinde paths. The first is
the ``multiplication'' map given by
\bea\label{P multiplication}
m_{\mathcal P}: {\mathcal P}_{k_1,l_1}^N\times{\mathcal P}_{k_2,l_2}^N&\to& {\mathcal P}_{k_1+k_2,l_1+l_2}^N,\notag\\
((\alb^{(1)};\beb^{(1)}),(\alb^{(2)};\beb^{(2)}))
&\mapsto&
(\alb^{(1)}+\alb^{(2)};\beb^{(1)}+\beb^{(2)}).
\ena
The map $m_{\mathcal P}$ is well-defined. This is obvious from the
definition of the admissibility (\ref{ADM}).

For an admissible triple $(l,l'',l')$, we have a well-defined injective left 
concatenation map
\bea\label{P concut}
c_{\mathcal P}(l,l'',l'): {\mathcal P}_{k,l'}^N\to {\mathcal P}_{k,l}^{N+1},\qquad
(\alb;{\beb})\mapsto (l,\alb;l'',{\beb}).
\ena

Note that multiplication and concatenation commute:
\bea\label{commute}
m_{\mathcal P}(c_{\mathcal P}(l_1,l_1'',l_1'),c_{\mathcal P}(l_2,l_2'',l_2'))=
  c_{\mathcal P}(l_1+l_2,l_1''+l_2'',l_1'+l_2') m_{\mathcal P}.
\ena

We have a recursive formula for Verlinde paths.
\bea\label{verlinde recursion}
{\mathcal P}_{k,l}^{N+1}=\bigsqcup_{l'',l':\,(l,l'',l')\;\;
{\rm is}\;\;{\rm admissible}} c_{\mathcal P}(l,l'',l'){\mathcal P}_{k,l'}^N.
\ena

\begin{lemma} 
The multiplication map $m_{\mathcal P}$ is surjective.
\end{lemma}
\begin{proof}
The statement is clear for $N=2$, i.e., for Verlinde triples. The general case follows from \Ref{commute} and \Ref{verlinde recursion}.
\end{proof}

\begin{lemma}\label{verlinde cardinality}
The cardinality of set ${\mathcal P}_{k,l}^N$ is given by the Verlinde number $d_{k,l}^{(N)}$ defined by \Ref{verlinde number}.
\end{lemma}
\begin{proof}
By \Ref{verlinde multiplication II} and \Ref{verlinde recursion}, the numbers computing the cardinality of ${\mathcal P}_{k,l}^N$ satisfy the same recursion relation as the Verlinde numbers $d_{k,l}^{(N)}$.
\end{proof}

\subsection{Recursion of combinatorial paths}
In this section we establish properties of the set of combinatorial
paths ${\mathcal C}^{(N)}_{k,l}$ similar to the properties of the set of
Verlinde paths ${\mathcal P}^{N}_{k,l}$ described in the previous
section.

We have an obvious ``multiplication'' map (compare to 
\Ref{P multiplication}):
\bea\label{C multiplication}
m_{\mathcal C}: {\mathcal C}_{k_1,l_1}^{(N)}\times{\mathcal
  C}_{k_2,l_2}^{N}&
\to& {\mathcal C}_{k_1+k_2,l_1+l_2}^{(N)},\notag\\
(({\bf a}^{(1)};{\bf b}^{(1)}),({\bf a}^{(2)};{\bf b}^{(2)}))&\mapsto&
({\bf a}^{(1)}+{\bf a}^{(2)};{\bf b}^{(1)}+{\bf b}^{(2)}).
\ena
We will prove in Theorem \ref{multiplication theorem} that it is surjective.

Next, we would like to obtain maps $c_{\mathcal C}(l,l'',l')$ similar
to $c_{\mathcal P}(l,l'',l')$ in \Ref{P concut}. We require three
properties for these maps. First, we would like to have a property
similar to the recursion relation \Ref{verlinde recursion}:
\bea\label{combinatorial recursion}
{\mathcal C}_{k,l}^{(N+1)}=\bigsqcup_{l'',l':\,(l,l'',l')\;\; {\rm is}\;\;
{\rm admissible}} c_{\mathcal C}(l,l'',l') {\mathcal C}_{k,l'}^{(N)}.
\ena
Second, we would like the maps $c_{\mathcal C}(l,l'',l')$ to be
extensions of the form
\bea\label{C concut}
c_{\mathcal C}(l',l'',l):
{\mathcal C}_{k,l'}^{(N)}\to {\mathcal C}_{k,l}^{(N+1)},\qquad ({\bf a};{\bf
  b})\mapsto ({\bf a},a;{\bf b},b),
\ena
where $a,b$ depend only on
the Verlinde triple $(l,l'',l')$ and $a_0$. This condition is natural
because it is similar to \Ref{P concut}.

The properties \Ref{combinatorial recursion}, \Ref{C concut}
uniquely determine the $c_{\mathcal C}(l,l'',l')$ for level $1$.
However, we need an additional property to fix the maps for the
general level.

We require the surjectivity of the
multiplication map. However, it does not uniquely fix the map.
We choose {\it a priori} the map (\ref{C concut})
so that Proposition \ref{multiplication theorem} holds.
This choice will be rederived {\it a posteriori} when we study
the filtration of Heisenberg modules in Section \ref{SEC-COIN}.

Now we describe such maps.

For each admissible triple of level $k$, $(l,l'',l')$,
we define maps $c_{\mathcal C}(l,l'',l')$ by formula \Ref{C concut} where
\bea\label{theta formula}
a&=&a(l,l'',l';a_0)=y(l,l'',l'),\label{RCOM1}\\
b&=&b(l,l'',l';a_0)=z(l,l'',l')-\Bigl(a_0-x(l,l'',l')\Bigr)^+.\label{RCOM2}
\ena
Here $t^+$ means the positive part of $t$: 
\bea\label{theta}
t^+=\max(t,0),
\ena
and $x,y,z$ are given by \Ref{xyz}.

\begin{proposition}\label{comb rec theorem}
The maps $c_{\mathcal C}(l,l'',l')$ given by the formula (\ref{C concut})
and (\ref{theta formula}) are well-defined. The combinatorial paths
enjoy the recursion (\ref{combinatorial recursion}).
\end{proposition}
\begin{proof}
For $({\bf a};{\bf b})\in{\mathcal C}_{k,l'}^{(N)}$ we have to check that
$({\bf a},a;{\bf b},b)\in{\mathcal C}_{k,l}^{(N+1)}$, where $a,b$ are given
by \Ref{theta formula}. First, $a=y(l,l'',l')\leq l$,
because $a=l-x(l,l'',l')$. Among the conditions \Ref{triangle condition}
for $({\bf a},a;{\bf b},b)$
we have new cases only for the triples $(a_0,b,a)$ and $(b_1,a_0,b)$.
They follow from
$$
a+b+a_0=y+z+a_0-(a_0-x)^+\leq y+z+x\leq k,
$$
$$
a_0+b_1+b=a_0+b_1+z-(a_0-x)^+\leq b_1+z+x=b_1+l'\leq k.
$$
(The last inequality follows from \Ref{trapezoid condition} with $i=-1,j=1$.)
Among the conditions \Ref{trapezoid condition} for $({\bf a},a;{\bf b},b)$
we have new cases only for $i=1,0,-1$. The case $i=0$ is weaker than $i=1$.
The case $i=1$ (resp., $i=-1$) follows from $b\leq l'$
(resp., $b+l\leq l'+a$) and the condition \Ref{trapezoid condition}
for $({\bf a};{\bf b})$ with $i=-1$.

Next, we have to check that the images of the maps $c_{\mathcal C}(l,l'',l')$ 
are disjoint. Suppose they are not, and
$a(l,l'',l';a_0)=a(l,l_1'',l'_1;a_0)$ and 
$b(l,l'',l',a_0)=b(l,l_1'',l'_1,a_0)$ for some $l,a_0$,$l'',l_1''$,$l',l'_1$. 
Then, the first equality implies $l''-l'=l_1''-l'_1$. Therefore, we have 
$x(l,l'',l')=x(l,l_1'',l_1')$. It follows that the positive part terms,
i.e., $(a_0-x)^+$, in the second
equality are equal and after its cancellation $l''+l'=l_1''+l'_1$. Therefore, 
we have $l''=l_1''$ and $l'=l_1'$.

Finally, let $({\bf a},a;{\bf b},b)\in{\mathcal C}_{k,l}^{(N+1)}$. Note that
$$
({\bf a},a,{\bf b},b)=c_{\mathcal C}(l,l'',l')({\bf a};{\bf b}),
$$
where
\bea
&&l'=x+z=l-a+b+(a_0-(l-a))^+,\notag\\
&&l''=y+z=a+b+(a_0-(l-a))^+.\notag
\ena
Moreover, the triple $(l,l'',l')$ is admissible, because 
$$
x+y+z=l+b+(a_0-(l-a))^+\leq a_0+a+b\leq k,
$$
by condition \Ref{triangle condition} for $({\bf a},a;{\bf b},b)$. Therefore,
it suffices to show that $({\bf a};{\bf b})\in{\mathcal C}_{k,l'}^{(N)}$.

This is also straightforward. The combinatorial path $({\bf a};{\bf b})$ looks
like
\bea
\left( \begin{array}{llllllllllllll}\notag
b_\infty=k&\cdots&0&&b_{N-1}&\cdots&&b_2&&b_1&&0&&l'\\
\dots&0&&a_{N-1}&&\cdots&a_2&&a_1&&a_0&&0
\end{array}\right),
\ena
and the combinatorial path $({\bf a},a;{\bf b},b)$ looks like 
\bea
\left( \begin{array}{llllllllllllllll}
b_\infty=k&\cdots&0&&b_{N-1}&\cdots&&b_2&&b_1&&b&&0&&l\\
\dots&0&&a_{N-1}&&\cdots&a_2&&a_1&&a_0&&a&&0
\end{array}\right).
\ena
The relation $a_0\leq l'$ follows from $b\geq 0$.
Conditions \Ref{triangle condition} are clearly preserved. Condition
\Ref{trapezoid condition} for $({\bf a};{\bf b})$ is also clear in all
places except when $i=0,-1$; moreover the case $i=0$ obviously follows
from the case $i=-1$. If $l\leq a+a_0$ then we have $l'=b+a_0$ and condition
\Ref{trapezoid condition} for $({\bf a};{\bf b})$ with $i=-1$ follows from
condition  \Ref{trapezoid condition} for $({\bf a},a;{\bf b},b)$ with $i=1$.
If $l\geq a+a_0$ then we have 
$l'+a=b+l$ and condition  \Ref{trapezoid condition} for $({\bf a};{\bf b})$
with $i=-1$ follows from condition  \Ref{trapezoid condition} for
$({\bf a},a;{\bf b},b)$ with $i=-1$. 
\end{proof}

Let ${\mathcal C}_{k,l}^{(N)}[i]\subset {\mathcal C}_{k,l}^{(N)}$
be the set of all combinatorial paths such that $a_0=i$. We have
\be
{\mathcal C}_{k,l}^{(N)}=\bigsqcup_{i=0}^l {\mathcal C}_{k,l}^{(N)}[i]. 
\ee
The recursion (\ref{combinatorial recursion}) splits into
\bea\label{comb recur}
{\mathcal C}_{k,l}^{(N+1)}[i]
=\bigsqcup_{l',l'',i'}c_{\mathcal C}(l,l'',l'){\mathcal C}_{k,l'}^{(N)}[i'],
\ena
where the sum is over $l',l'',i'$ such that $(l,l'',l')$ is Verlinde triple,
$i=(l+l''-l')/2$ and $i'\leq l'$.

Next we turn to the multiplication map. The analog of
the commutative property \Ref{commute} does not hold.
However we have the following
\begin{proposition}\label{multiplication theorem}
The multiplication map $m_{\mathcal C}$ given by \Ref{C multiplication}
is surjective.  Moreover,  let $({\bf a};{\bf b})\in{\mathcal C}_{k,l'}^{(N)}$
and let  $(l,l'',l')$ be an admissible triple of level $k$.
Let $(l,l'',l')=(l_1,l_1'',l'_1)+(l_2,l_2'',l'_2)$, where
$(l_1,l_1'',l'_1),(l_2,l_2'',l'_2)$ are admissible triples of levels
$k_1,k_2$, $k_1+k_2=k$. Then there exist
$({\bf a}^{(1)};{\bf b}^{(1)})\in{\mathcal C}_{k_1,l'_1}^{(N)}$,
$({\bf a}^{(2)};{\bf b}^{(2)})\in{\mathcal C}_{k_2,l'_2}^{(N)}$ such that 
$$
m_{\mathcal C}(({\bf a}^{(1)};{\bf b}^{(1)}),({\bf a}^{(2)};{\bf b}^{(2)}))
=({\bf a};{\bf b}),
$$
$$
m_{\mathcal C}(c_{\mathcal C}(l_1,l_1'',l'_1)
({\bf a}^{(1)};{\bf b}^{(1)}),c_{\mathcal C}(l_2,l_2'',l'_2)
({\bf a}^{(2)};{\bf b}^{(2)}))=c_{\mathcal C}(l,l'',l')({\bf a};{\bf b}).
$$ 
\end{proposition}
\begin{proof}
First, we prove that for any combinatorial path
$({\bf a};{\bf b})\in{\mathcal C}_{k,l'}^{N}$ and for any level $k$
Verlinde triple
$(l,l'',l')$, we have a decomposition of $({\bf a};{\bf b})$,
\bea
({\bf a};{\bf b})=\sum_{j=1}^k({\bf a}_j;{\bf b}_j)\notag,
\ena
where
$({\bf a}_j;{\bf b}_j)\in{\mathcal C}_{1,l'_j}^{(N)}$, $j=1,\dots,k$ and
$l_1'+\dots+l_k'=l'$; and a decomposition of the triple $(l,l'',l')$,
\bea
(l,l'',l')=\sum_{j=1}^k(l_j,l_j'',l'_j),\notag
\ena
where
$(l_j,l_j'',l'_j)$, $j=1,\dots,k$, are level $1$ triples, such that
we have
\bea\label{complete decomposition}
\sum_{j=1}^k  c_{\mathcal C}(l_j,l_j'',l'_j)({\bf a}_j;{\bf b}_j)
=c_{\mathcal C}(l,l'',l')({\bf a};{\bf b}).
\ena

We prove this statement by induction on $N$ using the induction hypothesis
that any level $k$ combinatorial path of length $N$ can be represented as
a sum of $k$ combinatorial paths of level $1$. It is obvious that
the decomposition (\ref{complete decomposition}) for $N$ implies
the above hypothesis for $N+1$.

The decomposition of $({\bf a};{\bf b})$ is not unique. We can choose any
decomposition $({\bf a}_j;{\bf b}_j)\in{\mathcal C}_{1,l'_j}^{(N)}$. 
Among the pairs $(l'_j,(a_j)_0)$, we have $a_0$ pairs $(1,1)$,
$l'-a_0$ pairs $(1,0)$ and $k-l'$ pairs $(0,0)$. We call these boundary pairs.
They are independent of the choice of the decomposition.

The decomposition of $(l,l'',l')$ is unique. We have $x$ triples $(1,0,1)$,
$y$ triples $(1,1,0)$, $z$ triples $(0,1,1)$ and $k-x-y-z$ triples $(0,0,0)$,
where $x,y,z$ are given by \Ref{xyz}.
We want to label these triple as $(l_j,l_j'',l'_j)$ so that the equality
\Ref{complete decomposition} holds. 

Recall \Ref{theta formula}.
The key property in the proof is that for level $1$
we have $a(1,1,0;0)=b(0,1,1;0)=1$ and $a=b=0$ in all other cases.

Now we match up the level $1$ triples with the boundary pairs. The triples
$(0,0,0)$ and $(1,1,0)$ must match up with the pairs $(0,0)$ because
of the consistency of $l'_j$. Their numbers are the same, and
we match up arbitrarily among them.

If $a_0\ge x$ we match up all $(1,0,1)$ with $(1,1)$, and we match up
$(0,1,1)$ with $(1,0)$ or $(1,1)$. If $a_0\le x$ we match up all $(0,1,1)$
with $(1,0)$, and we match up $(1,0,1)$ with $(1,0)$ or $(1,1)$.

It is easy to see \Ref{complete decomposition} by using the above remark
on $a,b$ for level $1$.

Now, we prove the statement of the proposition using the result obtained above.

One can group ${\mathcal L}=\{(l_j,l''_j,l'_j)\}_{j=1,\dots,k}$ as
${\mathcal L}={\mathcal L}^{(1)}\sqcup{\mathcal L}^{(2)}$ where
${\mathcal L}^{(i)}=\{(l^{(i)}_j,l^{(i)''}_j,l^{(i)'}_j)\}_{j=1,\dots,k_i}$
so that
$(l_i,l''_i,l'_i)=\sum_{j=1,\dots,k_i}(l^{(i)}_j,l^{(i)''}_j,l^{(i)'}_j)$.
We have the corresponding grouping
$\{({\bf a}^{(i)}_j,{\bf b}^{(i)}_j)\}_{j=1,\dots,k_i}$.
Then, setting
$({\bf a}^{(i)},{\bf b}^{(i)})=\sum_{j=1,\dots,k_i}
({\bf a}^{(i)}_j,{\bf b}^{(i)}_j)$, we have
$$
\sum_{j=1}^{k_i}c_{\mathcal C}(l^{(i)}_j,l^{(i)''}_j,l^{(i)'}_j)
({\bf a}^{(i)}_j;{\bf b}^{(i)}_j)
=c_{\mathcal C}(l_i,l''_i,l'_i)({\bf a}^{(i)};{\bf b}^{(i)}).
$$
The statement follows from this.
\end{proof}

\begin{remark}
As we noted, the surjectivity of the multiplication map
does not determine the map $c_{\mathcal C}(l,l'',l')$ uniquely.
On the other hand it is {\it not} possible to define a map
satisfying the commutativity similar to (\ref{commute}).
This requirement is too strong.
Proposition \ref{multiplication theorem} gives
the approriate requirement in order to fix the map
$c_{\mathcal C}(l,l'',l')$.
\end{remark}

\subsection{Bijection between the combinatorial paths and the Verlinde paths}
\label{bijection section}
In this section we use the results of the 
previous section to construct bijections $\iota$ between the sets of Verlinde
paths and the sets of combinatorial paths.

Define a map $\iota_{k,l}^N:{\mathcal P}_{k,l}^N\to {\mathcal C}_{k,l}^{(N)}$
as follows. We map $(\alb;\beb)$ to $({\bf a};{\bf b})$ with
\bea\label{bijection formula}
&&a_{i-1}=y(\al_i,\beta_i,\al_{i+1})
\qquad \qquad\qquad\qquad\qquad \qquad\qquad\qquad\;\;\;\;
(i=1,\dots,N),\notag\\ 
&&b_i=z(\al_i,\beta_i,\al_{i+1})
-(y(\al_{i+1},\beta_{i+1},\al_{i+2})-x(\al_i,\beta_i,\al_{i+1}))^+
\quad(i=1,\dots,N-1),
\ena
where $x,y,z$ are given by \Ref{xyz} and $(\,\,)^+$ by \Ref{theta}.

\begin{theorem}\label{bijection theorem}
The map $\iota_{k,l}^N:{\mathcal P}_{k,l}^N\to {\mathcal C}_{k,l}^{(N)}$ is
well-defined and bijective.
\end{theorem}
\begin{proof}
By the construction $\iota_{k,l}^N$ intertwines the maps 
$c_{\mathcal P}(l,l'',l')$ and $c_{\mathcal C}(l,l'',l')$.
Now the statement follows from
\Ref{verlinde recursion} and \Ref{combinatorial recursion} by induction on $N$.
\end{proof}

\begin{corollary}\label{CARLEM}
The cardinality of ${\mathcal C}_{k,l}^{(N)}$ is equal to the
Verlinde number $d_{k,l}^{(N)}$.
\end{corollary}
\begin{proof}
This follows from Theorem \ref{bijection theorem} and Lemma 
\ref{verlinde cardinality}.
\end{proof}

We define a triple grading on the Verlinde paths by using
the functions $e,s_1,s_2:{\mathcal P}_{k,l}^N\to \Z_{\geq 0}$ given by
\bea\label{energy}
e(\alb;\beb)&=&\sum_{i=1}^{N-1} i(a_i+b_i),\notag\\
s_1(\alb;\beb)&=&\sum_{i=1}^{N-1} b_i,\\
s_2(\alb;\beb)&=&a_0+\sum_{i=1}^{N-1} (a_i+b_i),\notag
\ena
where $a_i,b_i$ are given by \Ref{bijection formula}.

Note that the grading is local in the sense that
the grading functions are given by sums over terms
$a_i,b_i$ which depend only on two neighboring Verlinde triples
$(\alpha_i,\beta_i,\alpha_{i+1})$ and
$(\alpha_{i+1},\beta_{i+1},\alpha_{i+2})$.
 It is also intriguing that the form of $e$ is similar to the one-dimensional
configuration sum appearing in the corner transfer matrix method
in solvable lattice models (see, e.g., \cite{B}, \cite{DJKMO}).

Define the corresponding character
$$
\chi_{k,l}^{N}(q,z_1,z_2)= \sum
_{(\alb;\beb)\in{\mathcal P}_{k,l}^N}
q^{e(\alb;\beb)}z_1^{s_1(\alb;\beb)}z_2^{s_2(\alb;\beb)}.
$$ 
In Section \ref{SEC-COIN} we show that
$$
\chi_{k,l}^{(0,N)}(q,z)=z^l\chi_{k,l}^{N}(q,z^2,z^{-2}), 
$$
where
$\chi_{k,l}^{(0,N)}(q,z)$ is defined by \Ref{coinv character def}. The
meaning of the grading by $z_1,z_2$ will be clarified in Section
\ref{char identities section}.

\subsection{Recursion relations}\label{char rec section}
In this section we describe recursion relations for the characters.
We consider two types of the recursion relations which correspond to
the insertion of a Verlinde triple either at the beginning of the
Verlinde path (left concatenation) or at the end (right concatenation).

First, let us consider the right concatenation. Assume that $N\geq2$, and
consider $(\alb,\beb)\in{\mathcal P}_{k,l}^N$ and an admissible triple
$(\al,\beta,\gamma)$ such that $\al=\al_N$. Then the differences
of the gradings of the original path $(\alb;\beb)$ and the
concatenated path $(\alb,\gamma;\beb,\beta)$ depend only on the
new triple $(\al,\beta,\gamma)$ and the difference $\beta_{N-1}-\al_{N-1}$.
Note that $\al+\beta_{N-1}-\al_{N-1}=2z(\al_{N-1},\beta_{N-1},\al)$ as
in \Ref{xyz}.

Therefore, we split $\chi_{k,l}^{N}(q,z_1,z_2)$ into terms
which depend on the values of $\al_N$ and $z(\al_{N-1},\beta_{N-1},\al)$.
Namely, we define the partial character
$$
\chi_{k,l}^{N}[*;i,j](q,z_1,z_2)=
\sum _{\overset{\scriptstyle
(\alb;\beb)\in{\mathcal P}_{k,l}^{N+},}{{\al_{N}=i,\;\;
\al_N+\beta_{N-1}-\al_{N-1}=2j}}} q^{e(\alb;\beb)}z_1^{s_1(\alb;\beb)}
z_2^{s_2(\alb;\beb)}\quad(N\geq2).
$$

The partial character $\chi_{k,l}^{N}[*;i,j]$
is nontrivial only for $0\leq j\leq i\leq k$. 

\begin{proposition}
The partial characters $\chi_{k,l}^{N}[*;i,j](q,z_1,z_2)$
satisfy the following recursion relation
\be
\chi_{k,l}^{N+1} [*;i,j](q,z_1,z_2) =
\sum_{i',j'} R_k^{(N)}((i,j),(i',j'))\chi_{k,l}^{N}[*;i',j'](q,z_1,z_2),
\ee
where
$$
R_k^{(N)}((i,j),(i',j'))=q^{Ni-(N-1)(i-j-j')^+}z_1^{j'-(j+j'-i)^+}
z_2^{i-(i-j-j')^+},
$$
if $(i',2j-i+i',i)$ is admissible, and $R_k^{(N)}((i,j),(i',j'))=0$
otherwise.
\end{proposition}
\begin{proof}
Follows from \Ref{energy}.
\end{proof}

\begin{example}
For $k=1$, the matrix $R_1^{(N)}$ is given by
$$
\bordermatrix{&(0,0)&(1,0)&(1,1)\cr
(0,0)&1&1&1\cr
(1,0)&0&q&q^Nz_1z_2\cr
(1,1)&q^Nz_2&0&0\cr}.
$$
\end{example}

Similarly, for the left concatenation (see  \Ref{P concut})
we define the partial character
\bea\label{partial char}
\chi_{k,l}^{N}[i;*](q,z_1,z_2)=
\sum _{\overset{\scriptstyle
(\alb;\beb)\in{\mathcal P}_{k,l}^N,}{\al_1+\beta_1-\al_2=2i}}
 q^{e(\alb;\beb)}z_1^{s_1(\alb;\beb)}z_2^{s_2(\alb;\beb)}\quad(N\geq1).
\ena

The partial character $\chi_{k,l}^{N}[i;*]$ is nontrivial only
for $i\in\{0,\dots,l\}$.

\begin{proposition}
The partial characters $\chi_{k,l}^{N}[i;*](q,z_1,z_2)$ satisfy the following recursion relation
\bea\label{shift recursion}
\chi_{k,l}^{N+1} [i;*](q,z_1,z_2) =
\sum_{l',i'} L_k((l,i),(l',i'))\chi_{k,l'}^{N}[i';*](q,z_1,qz_2),
\ena
where
\be
L_k((l,i),(l',j))=(qz_1z_2)^{l'+i-l-(i'+i-l)^+}
q^{i'}z_2^i,
\ee
if $(l,2i-l+l',l')$ is admissible, and $L_k((l,i),(l',i'))=0$ otherwise.
\end{proposition}
\begin{proof}
Follows from \Ref{energy}.
\end{proof}

\begin{example}
For $k=1$, the matrix $L_1$ is given by
$$
\bordermatrix{&(0,0)&(1,0)&(1,1)\cr
(0,0)&1&qz_1z_2&q\cr
(1,0)&0&1&q\cr
(1,1)&z_2&0&0\cr}
$$
\end{example}

Note that the matrices $R_k^{(N)}$ and $L_k$ are nondegenerate square matrices
whose size is $k(k+1)/2$. The matrix $L_k$ does not depend on $N$, but the
corresponding recursion relation mixes different values of $l$ and includes
shifts with respect to $z_2$. The matrix $R_k^{(N)}$ does not mix $l$ but
depends on $N$.
Note also that at $q=z_1z_2=1$ their ranks degenerate to $k+1$. In this special
case our recursion relation is essentially the recursion relation for
Verlinde numbers $d_{k,l}^N$ which is given by multiplication by an element
in $k+1$ dimensional Verlinde algebra. This is called the Verlinde rule.
In this sense, our recursion relation can be thought of as a graded version of
of the Verlinde rule. 

Now we describe the relation between the two recursions above. The key
observation is the following lemma.
\begin{lemma}\label{inversion}
Let $(\alb;\beb)=(\alpha_1,\ldots,\alpha_N;\beta_1,\ldots,\beta_{N-1})$  be a
Verlinde path of length $N$.
Then 
$$
(\bar{\alb};\bar{\beb})=(\al_N,\al_{N-1},\dots,\al_1;\beta_{N-1},\dots,\beta_1)
$$
is also a Verlinde path of length $N$ and
\bea
e(\bar{\alb};\bar{\beb})+e(\alb;\beb)&=&(N-1)\,(s_2(\alb;\beb)+l-i),\notag\\
s_1(\bar{\alb};\bar{\beb})+l'-i'&=&s_1(\alb;\beb)+l-i,\notag\\
s_2(\bar{\alb};\bar{\beb})+l'-i'&=&s_2(\alb;\beb)+l-i,\notag
\ena
where 
$$
l'=\al_N,\qquad l=\al_1,\qquad i'=\half(\al_N+\beta_{N-1}-\al_{N-1}),\qquad i=\half(\al_1+\beta_1-\al_2).
$$
\end{lemma}
\begin{proof}
The assertion follows from formula \Ref{energy}.
Indeed,
$$
a_i+b_i=\half(\beta_i+\beta_{i+1}+2\al_{i+1}-\al_i-\al_{i+2}
-(\beta_i+\beta_{i+1}-\al_i-\al_{i+2})^+) =\bar{a}_{N-1-i}+\bar{b}_{N-1-i},
$$
$i=1,\dots,N-2$. Also, we always have $b_0=b_{N-1}=0$. Therefore, for energy $e$ and spin $s_2$ we only have to compensate the contribution from $a_{N-1}$ and $a_0$.

To prove the equality for $s_1$ we note that $$
s_2(\alb,\beb)-s_1(\alb,\beb)=\half(\al_1+\al_N+\sum_{j=1}^{N-1} \beta_j)
=s_2(\bar{\alb},\bar{\beb})-s_1(\bar{\alb},\bar{\beb}).
$$
\end{proof}

Note that for the Verlinde path $(\alb,\beb)$ in Lemma
\ref{inversion}, we have $l=\bar a_{N-1}$, $l'=a_{N-1}$, $i=a_0$ and
$i'=\bar a_0$.

We decompose the partial characters further and define
$$
\chi_{k}^{N}[i,l;i',l'](q,z_1,z_2)=
\sum_{\overset{\scriptstyle
(\alb;\beb)\in{\mathcal P}_{k,l}^N,\;\;\al_1+\beta_1-\al_2=2i}{
  \al_N=l',\;\; \al_N+\beta_{N-1}-\al_{N-1}=2i'}}
q^{e(\alb;\beb)}z_1^{s_1(\alb;\beb)}z_2^{s_2(\alb;\beb)}\quad(N\geq2).
$$
Then, we have 
$$
\sum_{i=0}^l
\chi_{k}^{N}[i,l;i',l']=\chi_{k,l}^{N}[*;i',l'], 
$$
$$
\sum_{\overset{\scriptstyle i',l'=0}{ i'\leq l'}}^k
\chi_{k}^{N}[i,l;i',l']=\chi_{k,l}^{N}[i;*], 
$$ 
$$
\sum_{\overset{\scriptstyle i',l',i=0}{ i'\leq l',\;\; i\leq l}}^k
\chi_{k}^{N}[i,l;i',l'] =\chi_{k,l}^N.  
$$
Then, left and right concatenations act on vectors of the form
$$
\{\chi_{k}^{N}[i,l;i',l']\}
_{i,l,i',l'\in\{1,\dots,k\},\;0\leq i\leq l,\;0\leq i'\leq l'} 
$$ 
as block diagonal matrices with $k(k+1)/2$ blocks,
each being equal to $L_{k}$ and $R_k^{(N)}$ (plus a shift $z_2\mapsto
qz_2$) respectively. By Lemma \ref{inversion}, these two actions are
conjugated by a transformation
$\chi_{k}^{N}[i,l;i',l']\to\chi_{k}^{N}[i',l';i,l]$, 
$$
\chi_{k}^N[i,l;i',l'](q,z_1,z_2)=q^{(N-1)(l'-i')}(z_1z_2)^{l'-i'-l+i}\,
\chi_{k}^N[i',l';i,l](q^{-1},z_1,z_2q^{N-1}).  
$$

In the subsequent paper we will present
explicit fermionic formulas for the characters in terms of $q$-binomial
coefficients.

\newcommand{\bU}{\overline{U}}
\newcommand{\bV}{\overline{V}}
\newcommand{\bW}{\overline{W}}
\section{Heisenberg modules}\label{SEC3}
In this section we discuss the reduction of the problem of computing
coinvariants for $\wsl$ to the problem of computing coinvariants for
affine Heisenberg algebra ${\widehat{\mathcal{H}}}$.
Then we begin a study of the structure of "integrable" Heisenberg modules.

\subsection{Abelianization}\label{AHA} 
In this section we discuss the idea of passage from $\wsl$ module
to $\cH$ module.

Let us recall that the Verma module $M_{k,l}$ is a quotient of the free
module generated by $v$ over the universal enveloping algebra of $\wsl$ by  
the highest weight conditions for $v$.
The highest weight conditions in the Verma module $M_{k,l}$ have the form
$e_iv=0$ $(i\le0)$, $h_iv=0$ $(i\le-1)$, $h_0v=lv$
and $f_iv=0$ $(i\le-1)$. The integrable module $L_{k,l}$ is the quotient of
the Verma module $M_{k,l}$ by the submodule generated by the elements
of the form  $e(z)^{k+1}w,f(z)^{k+1}w$ for $w\in M_{k,l}$.
The space of $(M,N)$-coinvariants
$L^{(M,N)}_{k,l}$ is the further quotient of $L_{k,l}$ by ${\rm Im}\,e_i$
$(i\geq M)$ and ${\rm Im}\,f_i$ $(i\geq N)$. 

We call a vector in $L^{(M,N)}_{k,l}$ of the form $\m v$ where $\m$
is a monomial in $U(\wsl)$ as a monomial vector or simply as a monomial.
Our goal is to obtain a monomial basis of $L^{(M,N)}_{k,l}$, i.e.,
a set of monomial vectors in $L^{(M,N)}_{k,l}$ which forms a basis.
In principle, one can rewrite given monomials to linear
combinations of other monomials by using the above relations in the hope that
we can find a monomial basis as a result.
However, this is not a unique procedure, and in general, we do not know
any canonical way of doing such reduction, which eventually leads
to a monomial basis.

In \cite{FKLMM} an approach using a filtration of $L_{k,l}$
and its adjoint graded space ${\rm Gr}(L_{k,l})$
was used in the special case $k=1$. In this paper we further develop this
method.

Let us explain the basic idea in this approach.
If we consider a filtration such that each generators $e_i,f_i,h_i$
($i\in\Z$) of $\wsl$ has a fixed degree with respect to the filtration,
we can induce its graded action on ${\rm Gr}(L_{k,l})$. We consider the space
of coinvariants in ${\rm Gr}(L_{k,l})$ with respect to the induced actions
of $e_i$ ($i\geq M$) and $f_i$ ($i\geq N$). For an appropriate choice
of filtration, the relations in terms of the induced actions simplify,
and we may find a monomial basis of the space of coinvariants
${\rm Gr}(L_{k,l})^{(M,N)}$ by exploiting the induced actions.
If, in addition, the equality of dimensions
$$
{\rm dim}\,L_{k,l}^{(M,N)}={\rm dim}\,{\rm Gr}(L_{k,l})^{(M,N)}
$$
holds, the same set of monomials gives a basis of  $L_{k,l}^{(M,N)}$.

Now, let us discuss a concrete example of filtration.

Let $\mathfrak{g}$ be a Lie algebra with a Lie subalgebra
$\mathfrak{a}\subset \mathfrak{g}$ and $V$ a cyclic $\mathfrak{g}$-module
with a cyclic vector $v$. Consider the filtration
$A(V)=(A_i(V))_{i\geq0}$ defined inductively by 
\bean A_0(V)=\C v,\quad A_{i+1}(V)=A_i(V)+{\mathfrak{g}}\cdot A_i(V).  
\eean 
Then the induced actions of elements of  $\mathfrak{g}$ on ${\rm
Gr}^A(V)$ with respect to this filtration commute. We call the passage from
the module $V$ to ${\rm Gr}^A(V)$ the abelianization procedure.
Now by Lemma \ref{deformation} we have the inequality of dimensions of
coinvariants 
\bea 
{\rm
  dim}\,V/{\mathfrak{a}}\,V\leq {\rm dim}\,{\rm
  Gr}^A(V)/\overline{\mathfrak{a}}\,{\rm Gr}^A(V),
\label{INEQU}
\ena 
where $\overline{\mathfrak{a}}$ is the induced action of
$\mathfrak{a}$. 

In the case we consider in this paper, 
${\mathfrak{g}}=\widehat{\mathfrak{sl}}_2$,
${\mathfrak{a}}=\left(\oplus_{i\geq M}e_i\right)
\oplus\left(\oplus_{i\geq M+N}h_i\right)
\oplus\left(\oplus_{i\geq N}f_i\right)$ and $V=L_{k,l}$,
the equality in (\ref{INEQU}) does not hold. In other words, the
abelianization loses the information about the character of the original
coinvariants. The dimension of the coinvariants jumps in the abelian limit.

We will successfully exploit a finer filtration given
by the actions of $\{e_i\}$ and $\{f_i\}$ (but not $\{h_i\}$). The
associated graded spaces are Heisenberg modules, which are more
tractable than $\wsl$-modules. Moreover, it turns out that
the dimension of the coinvariants does not jump in this limit if $M,N\geq 1$. 
(In the cases $M=0,N>0$ and $M>0, N=0$ the dimension does not jump either if one chooses an appropriate cyclic vector in $L_{k,l}$.)

\subsection{Affine Heisenberg algebra}
In this section, we derive an action of the affine Heisenberg algebra
${\widehat{\mathcal{H}}}$ on the associated graded space of the
integrable $\wsl$-modules, filtered by the action of the currents
$e(z)$ and $f(z)$.

The affine Heisenberg algebra $\cH$ is a Lie algebra generated by the elements
$e_i,h_i,f_i$ $(i\in\Z)$ satisfying the relations
\bean
[e_i,f_i]\,=\,h_{i+j},\qquad [h_i,e_j]\,=\,[h_i,f_j]\,=\,0.
\eean

Here and below we use the same symbols $e_i$, $f_i$ and $h_i$ for the
generators of $\cH$ as those of $\wsl$.

An $\cH$-module $V$ is called a highest weight module if
it is generated by a highest weight vector $v$ satisfying
\bean
e_iv=0,\qquad h_iv=0,\qquad f_iv=0\qquad (i\le-1).
\eean
A highest weight $\cH$- module $V$ is called level-$k$ integrable if
the action satisfies
\bea
e(z)^{k+1}w=f(z)^{k+1}w=0\label{INT},
\ena
for any vector $w\in V$. 
For simplicity, we call a level $k$ integrable highest weight module
an $\cH_k$-module. 

Let $V$ be an $\cH_k$-module.
Since $({\rm ad}\,e_0)f(z)=h(z)$ and $({\rm ad}\,f_0)e(z)=-h(z)$, we
have the relations
\bea
e(z)^{k+1-a}h(z)^aw=f(z)^{k+1-a}h(z)^aw=0\quad(0\leq a\leq k+1)
\label{INT2}
\ena
for all $w\in V$.

The $\cH_k$ module $V$ is graded by
\bean
{\rm deg}\, v=(0,0,0),\quad {\rm deg}\,e_i=(1,0,i),\quad 
{\rm deg}\,f_i=(0,1,i),\quad
{\rm deg}\,h_i=(1,1,i),
\eean
where $v$ is the highest weight vector in $V$.
Let $V_{m,n,e}$ be the subspace of degree $(m,n,e)$ and $V_{m,n}=\oplus_eV_{m,n,e}$.

Below, we will introduce various $\cH$-linear mappings $\mu:V\rightarrow W$
among $\cH_k$-modules. All these maps are graded maps. If
the highest weight vector $v$ of $V$ is mapped to a vector
in a subspace of some degree $(m,n,e)$ then the mapping $\mu$ has
degree $(m,n,e)$.

Let $L_{k,l}$ be the integrable $\wsl$ module of level $k$ and highest
weight $l$. Define a filtration $F(L_{k,l})=(F_i(L_{k,l}))_{i\geq0}$
of $L_{k,l}$ inductively by
\bean
F_0(L_{k,l})=\C v_k[l],\quad F_{i+1}(L_{k,l})=F_i(L_{k,l})+
{\mathfrak{c}}\cdot F_i(L_{k,l}),
\eean
where
\bean
{\mathfrak{c}}&=&
\Bigl(\bigoplus\limits_{i\in\Z}\C e_i\Bigr)
\bigoplus\Bigl(\bigoplus\limits_{i\in\Z}\C f_i\Bigr)
\subset\wsl.
\eean
We define the induced action of $e_i$, $f_i$ and $h_i$ on the
associated graded space to have degree $1$ for $e_i$ and $f_i$, and $2$
for $h_i$.

\begin{proposition}\label{IND}
The induced actions of $e_i$, $f_i$, $h_i$ $(i\in\Z)$, define an
$\cH_k$-module structure on
the vector space ${\rm Gr}^F(L_{k,l})$.
\end{proposition}
\begin{proof}
Let us denote the induced actions by $\bar e_i$, $\bar f_i$, $\bar h_i$.
We will show $[\bar h_i,\bar e_j]=0$. (The proof of
$[\bar h_i,\bar f_j]=0$ is similar, and the rest of the relations are
straightforward.) We have $[h_i,e_j]=2e_{i+j}$. The mapping
$[\bar h_i,\bar e_j]$ sends ${\rm Gr}^F_n(L_{k,l})$ to
${\rm Gr}^F_{n+3}(L_{k,l})$. For $v\in F_n(L_{k,l})$, let $\bar v^{(n)}$
be the corresponding element in ${\rm Gr}^F_n(L_{k,l})$. Then we have
\be
[\bar h_i,\bar e_j]\bar v^{(n)}=\overline{[h_i,e_j]v}^{(n+3)}=2\overline{e_{i+j}v}^{(n+3)}=0.
\ee
\end{proof}

We define an involution $\iota$ of $\wsl$ by
\bea
\iota(e_i)=f_i,\qquad\iota(f_i)=e_i,\qquad\iota(h_i)=-h_i.
\label{INVOL0}
\ena 
and an involution of $L_{k,l}$ which we also denote $\iota$ by 
\bean
\iota(v_k[l])=u_k[l],\\
\iota(x\cdot y)=\iota(x)\cdot\iota(y),
\eean
for $x\in\wsl$ and $y\in L_{k,l}$. Here we set
\begin{equation}
u_k[l]=\frac{f_0^l}{l!}v_k[l].\label{UVEC}
\end{equation}

The involution $\iota$ induces a new filtration $E$ of $L_{k,l}$ from
$F$. Namely, it is given by a procedure similar to $F$ above, 
\be
E_0(L_{k,l})=\C u_k[l],\quad E_{i+1}(L_{k,l})=E_i(L_{k,l})+
{\mathfrak{c}}\cdot E_i(L_{k,l}).
\ee
The adjoint graded space ${\rm Gr}^E(L_{k,l})$ is also an $\cH_k$-module.
\subsection{Heisenberg modules}
In this section we define three types of $\cH_k$-modules which differ
only by the annihilation conditions for the
highest weight vector.

Let $l_1,l_2,l_3\geq-1$.  Let $V_k[l_1,l_2,l_3]$ be the $\cH_k$-module
generated by the highest weight
vector $v_V=v_V[l_1,l_2,l_3]$, satisfying the relations
\bea
&&e_iv_V=0\quad(i\leq0),\qquad f_iv_V=0\quad(i\leq-1),\label{HIGH}\\
&&e_1^{l_1+1}v_V=f_0^{l_2+1}v_V=h_0^{l_3+1}v_V=0.\label{REL1}
\ena

Let $U_k[l_1,l_2,l_3]$ be the $\cH_k$-module
generated by the highest weight
vector $v_U=v_U[l_1,l_2,l_3]$, satisfying the relations
\bea
&&e_iv_U=0\qquad(i\le-1),\qquad f_iv_U=0\qquad(i\leq0),\label{REL4}\\
&&e_0^{l_1+1}v_U=f_1^{l_2+1}v_U=h_0^{l_3+1}v_U=0.\label{REL5}
\ena

Let $W_k[l_1,l_2,l_3]$ be the $\cH_k$-module
generated by the highest weight
vector $v_W=v_W[l_1,l_2,l_3]$, satisfying the relations
\bea
&&e_iv_W=f_iv_W=0\qquad (i\leq0),\label{REL2}\\
&&e_1^{l_1+1}v_W=f_1^{l_2+1}v_W=h_1^{l_3+1}v_W=0.\label{REL3}
\ena

The modules $V_k[l_1,l_2,l_3]$, $U_k[l_1,l_2,l_3]$ and
$W_k[l_1,l_2,l_3]$ are related through simple automorphisms of $\cH$
as follows.

Define the shift automorphism $T$ of $\cH$ by
\bea
T(e_i)=e_i,\qquad T(f_i)=f_{i+1},\qquad T(h_i)=h_{i+1}.\label{SHIFT}
\ena  
We have an isomorphism of linear spaces
\bea\label{V=W}
T: V_k[l_1,l_2,l_3]\stackrel{\sim}{\to} W_k[l_1,l_2,l_3],
\ena
determined by the properties
\bean
T(v_V)&=&v_W,\\
T(x\cdot y)&=&T(x)\cdot T(y)\quad(x\in \cH,\ y\in V_k[l_1,l_2,l_3]).
\eean

Define involution $\iota$ of $\cH$ (cf. \Ref{INVOL0}) by
\bea
\iota(e_i)=f_i,\qquad \iota(f_i)=e_i,\qquad \iota(h_i)=-h_i.
\label{INVOL}
\ena
We have an isomorphism of linear spaces
\bea\label{V=U}
\iota: V_k[l_1,l_2,l_3] \stackrel{\sim}{\to} U_k[l_2,l_1,l_3], 
\ena
determined by the properties
\bean
\iota(v_V)&=&v_U,\\
\iota(x\cdot y)&=&\iota(x)\cdot\iota(y)\quad(x\in \cH,\ y\in V_k[l_1,l_2,l_3]).
\eean

Note that if $l_i=-1$ for some $i$, then we have 
\be
V_k[l_1,l_2,l_3]=U_k[l_1,l_2,l_3]=W_k[l_1,l_2,l_3]=0.
\ee

The formulas \Ref{INT}, \Ref{INT2} imply
\be
e_1^{k+1}v=f_0^{k+1}v=h_0^{k+1}v=0,
\ee
for any highest weight vector in a $\cH_k$-module.
Therefore, we have 
\be
V_k[l_1,l_2,l_3]=V_k[\min(l_1,k),\min(l_2,k),\min(l_3,k)].
\ee

{}From Lemma \ref{LM0} we obtain another relation
\be
V_k[l_1,l_2,l_3]=V_k[l_1,l_2;{\rm min}(l_1,l_2,l_3)].
\ee

Similarly we have
\bean
U_k[l_1,l_2,l_3]&=&U_k[\min(l_1,k),\min(l_2,k),\min(l_3,k)],\\
U_k[l_1,l_2,l_3]&=&U_k[l_1,l_2;{\rm min}(l_1,l_2,l_3)],
\eean
and
\bean
W_k[l_1,l_2,l_3]&=&W_k[\min(l_1,k),\min(l_2,k),\min(l_3,k)],\\
W_k[l_1,l_2,l_3]&=&W_k[l_1,l_2;{\rm min}(l_1,l_2,l_3)].
\eean

Introduce the following simplified notation:
\bean
\bV_k[l_1,l_2]&=&V_k[l_1,l_2;\min(l_1,l_2)],\\
\bU_k[l_1,l_2]&=&U_k[l_1,l_2;\min(l_1,l_2)],\\
\bW_k[l_1,l_2]&=&W_k[l_1,l_2;\min(l_1,l_2)].
\eean
These are the modules where $h$ condition on the highest weight vector in
\Ref{REL1},\Ref{REL5},\Ref{REL3} follows from
$e$ and $f$ conditions. So, these are the ``simplest'' and the largest modules.

Also introduce the following notations.
\bean
U_k[l_1,l_2]&=&U_k[l_1,l_2,0],\\
V_k[l_1,l_2]&=&V_k[l_1,l_2,0],\\
W_k[l_1,l_2]&=&W_k[l_1,l_2,l_1+l_2-k].
\eean
These are the modules which are directly related to the $\wsl$ modules
studied in this paper (see Section \ref{connection section}).

\subsection{Filtration and recurrence relations}\label{FILTR}
In this section we will determine the structure of the $\cH_k$-module
$W_k[l_1,l_2,l_3]$ by recurrence relations. 
We construct a filtration of $W_k[l_1,l_2,l_3]$ such that
each graded component of the adjoint graded space is isomorphic to
$T(W_k[l_1',l_2'])$ for some $l_1',l_2'$ such that $0\le l_1',l_2'\le k$ and
$l_1'+l_2'\ge k$. 

There is a significant difference between the $\wsl$-modules
$L_{k,l}$ and the $\cH_k$-modules $U_k[l_1,l_2,l_3]$, $V_k[l_1,l_2,l_3]$
and $W_k[l_1,l_2,l_3]$. The
modules $L_{k,l}$ are irreducible, while the Heisenberg 
modules contain infinitely many submodules.
We are interested in the self-similarity of subquotients of
$\cH_k$-modules, which will allow us to construct a basis of
$\cH_k$-modules by a recursion procedure. 

First we consider a relatively simple case, when a submodule of a $\cH_k$ is isomorphic to 
$W_k[l_1',l_2']$ for some $l_1',l_2'$.

Let $v$ be the highest weight vector in a highest weight $\cH_k$-module.
We denote the submodule generated by a set of vectors $\m_iv$ $(i\in I)$
by $\langle \m_i\ (i\in I)\rangle$.

We have the following short exact sequences:
\bean
&&0\rightarrow \la f_0^{l_2}\ra\rightarrow
V_k[l_1,l_2]\rightarrow V_k[l_1,l_2-1]\rightarrow0\quad(l_1+l_2\leq k),
\label{54}\\
&&0\rightarrow \la e_0^{l_1}\ra\rightarrow
U_k[l_1,l_2]\rightarrow U_k[l_1-1,l_2]\rightarrow0\quad(l_1+l_2\leq k).
\label{540}
\eean
\begin{proposition}\label{triple prop}
The following are well-defined $\cH$-linear mappings.
\bea
W_k[l_1+l_2,k-l_2]&\rightarrow&V_k[l_1,l_2],
\quad \m v_W\mapsto \m f_0^{l_2}v_V\quad(l_1+l_2\leq k),\label{SH3}
\\
W_k[k-l_1,l_1+l_2]&\rightarrow&U_k[l_1,l_2],
\quad \m v_W\mapsto \m e_0^{l_1}v_U\quad(l_1+l_2\leq k),\label{SH30}
\ena 
\end{proposition}
\begin{proof}
We give the proof for the first map. The second one is similar.

We need to show the validity of the relations (\ref{REL2}) and (\ref{REL3}).
The relation $e_0f_0^{l_2}v=0$ follows from
$h_0v=0$. The relation $e_1^{l_1+l_2+1}f_0^{l_2}v=0$ follows from 
a variant of Lemma \ref{LM1}. Finally, the relation $f_1^{k-l_2+1}f_0^{l_2}v=0$
follows from the integrability condition (\ref{INT}), and
$h_1^{l_1+1}f_0^{l_2}v=0$ from a variant of Lemma \ref{LM3}.
\end{proof}

\begin{corollary}\label{TRIPLES}
The sequences
\bea
&&W_k[l_1+l_2,k-l_2]\rightarrow V_k[l_1,l_2]\rightarrow
V_k[l_1,l_2-1]\rightarrow0\quad(l_1+l_2\leq k)\label{triple2},\\
&&W_k[k-l_1,l_1+l_2]\rightarrow U_k[l_1,l_2]\rightarrow
U_k[l_1-1,l_2]\rightarrow0\quad(l_1+l_2\leq k)\label{triple5}
\ena
are exact.
\end{corollary}
\begin{proof} Corollary \ref{TRIPLES} follows from Proposition
\ref{triple prop}. 
\end{proof} 

Later, we will prove that \Ref{SH3} and \Ref{SH30} are in fact injective
(see Corollary \ref{SES cor2}).

Now we consider a more complicated situation.
The module $W_k[l_1,l_2,l_3]$ has a filtration, such that each of
the corresponding subquotients is isomorphic to an $\cH_k$-module of the form
$T(W_k[l_1',l_2'])$ for some $l_1',l_2'$.

Fix $l_1,l_2,l_3$ such that $0\le l_1,l_2\le k$ and
$0\leq l_3\leq{\rm min}(l_1,l_2)$, and denote the highest weight vector
of $W_k[l_1,l_2,l_3]$ by $v$.

\begin{lemma}
The vector $h_1^af_1^cv\in W_k[l_1,l_2,l_3]$ is zero unless
\bea
0\leq a\leq l_3,\quad0\leq c\leq l_2-a.
\label{REGION}
\ena
\end{lemma}
\begin{proof}
The relation $h_1^af_1^cv=0$ for $a>l_3$ follows from $h_1^{l_3+1}v=0$.
The relation $h_1^af_1^cv=0$ for $a+c>l_2$ follows from 
$f_1^{l_2+1}v=e_0v=0$ by a variant of Lemma \ref{LM0}.
\end{proof}

Now, we define a filtration of $W_k[l_1,l_2,l_3]$.
Consider submodules of the form $\la h_1^af_1^c\ra$.
Since these submodules are not linearly ordered with respect to inclusion,
we define a two-step filtration. First, define the filtration by
$\la h_1^a\ra$, then refine this filtration by using $h_1^af_1^c$.

Set
\bean
W_{a,c}&=&\langle h_1^af_1^c,h_1^{a+1}\rangle\subset W_k[l_1,l_2,l_3].
\eean 
Note that
\bean
W_{a,c}&\supset&W_{a,c+1},\\
W_{a,0}&=&\langle h_1^a\rangle,\\
W_{a,c}&=&W_{a+1,0}\hbox{ if $c>l_2-a$}.
\eean
Thus we have a filtration of $W[l_1,l_2,l_3]$:
\bea
W[l_1,l_2,l_3]=W_{0,0}\supset W_{0,1}\supset\cdots
\supset W_{1,0}\supset W_{1,1}\supset
\cdots\supset W_{l_3,0}\supset\cdots\supset W_{l_3,l_2-l_3+1}=0.
\label{FILT2}
\ena
We call the filtration of $W_k[l_1,l_2,l_3]$ given by (\ref{FILT2}) the
canonical filtration of the first kind, and denote it by $C(W_k[l_1,l_2,l_3])$.

\begin{proposition}\label{DECOM}
For each $(a,c)$ satisfying $(\ref{REGION})$ there exists an $\cH$-linear
surjection
\begin{eqnarray}
T(W_k[l_1',l_2'])&\rightarrow&{\rm Gr}^C_{a,c}(W_k[l_1,l_2,l_3]),\label{MAP}\\
T(\m v_W[l'_1,l'_2,l'_1+l'_2-k])&\mapsto&T(\m)v_{a,c},\nonumber\\
v_{a,c}&=&h_1^af_1^cv_W[l_1,l_2,l_3].\nonumber
\end{eqnarray}
Here
\bea\label{l primes}
l_1'={\rm min}(k-a,l_1+c-a),\quad l_2'=k-c\label{L'}.
\ena
\end{proposition}
\begin{proof}
It is enough to show that the vector $v_{a,c}$ satisfies
\bean
&&e_iv_{a,c}=0\quad (i\leq0),\qquad f_iv_{a,c}=0\quad(i\leq1),\\
&&e_1^{l_1'+1}v_{a,c}=0,\quad f_2^{l_2'+1}v_{a,c}=0,\quad
h_2^{l_1'+l_2'-k+1}v_{a,c}=0.
\eean
First, note that the quotient module ${\rm Gr}^C_{a,c}(W_k[l_1,l_2,l_3])$
is generated by $h_1^af_1^cv$ where $v=v_W[l_1,l_2,l_3]$. Note also that
that $h_1^af_1^{c+1}v=0$ and $h_1^{a+1}f_1^{c-1}v=0$.

If $i\le0$, then $e_iv=0$. Therefore, we have
\be
e_iv_{a,c}=e_ih_1^af_1^cv=ch_1^ah_{i+1}f_1^{c-1}v=0
\ee
because 
\bean
h_{i+1}v&=&0\qquad(i<0),\\
h_1^{a+1}f_1^{c-1}v&=&0\qquad (i=0).
\eean

If $i\le1$, we have
\be
f_iv_{a,c}=f_ih_1^af_1^cv=0\
\ee
because 
\bean
f_iv&=&0\qquad i\le0,\\
h_1^af_1^{c+1}v&=&0\qquad (i=1).
\eean

We have 
\bean
f_2^{l_2'+1}v_{a,c}=h_1^af_2^{k-c+1}f_1^cv=0
\eean
because of integrability condition (\ref{INT}).

Consider the case $0\le c\le k-l_1$. We have $l_1'=l_1+c-a$. From
Lemma \ref{LM2} we have
\bean
e_1^{l_1'+1}v_{a,c}=e_1^{l_1+c-a+1}h_1^af_1^cv=0.
\eean
{}From Lemma \ref{LM3} we have
\bean
h_2^{l_1'+l_2'-k+1}v_{a,c}=h_2^{l_1-a+1}h_1^af_1^cv=0.
\eean

Finally, consider the case $k-l_1\le c\le l_2-a$.
We have $l_1'=k-a$, and by using (\ref{INT2}), we obtain
\bean
&&e_1^{l_1'+1}v_{a,c}=e_1^{k-a+1}h_1^af_1^cv=0,\\
&&h_2^{l_1'+l_2'-k+1}v_{a,c}=h_2^{k-a-c+1}h_1^af_1^cv=0.
\eean
\end{proof}

In Section \ref{general recur section}, we will prove that the map (\ref{MAP}) is an
isomorphism. 

\section{Recursion for the coinvariants}\label{SEC-COIN}
In this section we prove the main theorem (Theorem \ref{MAINTH}),
which states that the dimension of the space of the $\wsl$ coinvariants
$L^{(M,N)}_{k,l}$ is given by the Verlinde rule.
The key idea is to derive a recursion relation for the spaces of the
Heisenberg coinvariants with respect to the parameters $M,N$.

\subsection{The spaces of $(M,N)$-coinvariants}
In Section \ref{defs} we have introduced the space of $(M,N)$-coinvariants
$L^{(M,N)}_{k,l}$. In this section we define similar spaces for
$\cH_k$-modules and prove that the properties of $\cH_k$ modules
described in Section \ref{FILTR} are valid for the coinvariants as well.

For $M,N\geq0$, we set
\bean
{\mathfrak{a}}^{(M,N)}&=&\left(\bigoplus\limits_{i\geq M}\C e_i\right)
\oplus\left(\bigoplus\limits_{i\geq M+N}\C h_i\right)
\oplus\left(\bigoplus\limits_{i\geq N}\C f_i\right)\subset \cH.
\eean
For an $\cH_k$-module $V$, the space of $(M,N)$-coinvariants is defined by
\bean
V^{(M,N)}&=&V/{\mathfrak{a}}^{(M,N)}V.
\eean

\begin{remark}
For simplicity of notation, we denote the image of the highest weight vector
$v$ of an $\cH_k$-module $V$ in the space of $(M,N)$-coinvariants by $v$.
\end{remark}

Define
\bea
W_k^{(M,N)}[l_1,l_2,l_3]&=&\bigoplus\limits_{m,n}W_k^{(M,N)}[l_1,l_2,l_3]_{m,n},\no\\
W_k^{(M,N)}[l_1,l_2,l_3]_{m,n}&=&
\begin{cases}
0\text{ if $M=0$ and $n-2m<k-l_1$};\\
0\text{ if $N=0$ and $m-2n<k-l_2$};\\
\oplus_e W_k[l_1,l_2,l_3]^{(M,N)}_{m,n,e}\text{ otherwise}.
\end{cases}\label{MD2}
\ena

Namely, we cut off some weight spaces of
$W_k[l_1,l_2,l_3]^{(M,N)}$ when either $M$ or $N$ is zero. 
Note, in particular, that
\bea
W_k^{(0,0)}[l_1,l_2,l_3]&=&
\begin{cases}
\C\cdot v\text{ if $l_1=l_2=k$;}\\
0\text{ otherwise.}
\end{cases}
\label{INIT}
\ena
Note also the difference in notation between
$W_k^{(M,N)}[l_1,l_2,l_3]$ and $W_k[l_1,l_2,l_3]^{(M,N)}$. They are
equal if $M,N>0$.

\begin{remark}
The reason for the above cutoff is to obtain the spaces satisfying
the recursion relations (\ref{MAP2-2}). One can easily check that
the cutoff condition is equivalent to the restriction
$\sum_ia_i-\sum_ib_i\geq l$ for the monomial basis (\ref{FHBAS})
following from (\ref{trapezoid condition}).
\end{remark}

Now we establish the coinvariant version of Corollary \ref{TRIPLES}.
First we deal with the cutoffs.

\begin{lemma}\label{LEM}
Consider the mapping
\bean
W_k[k-l_1,l_1+l_2]^{(0,N)}_{m,n,e}\rightarrow U_k[l_1,l_2]^{(0,N)}_{m+l_1,n,e},
\quad
{\bf m}v_W\mapsto {\bf m}e_0^{l_1}v_U
\eean
induced from (\ref{SH30}).
If $n-2m<l_1$, the image of this mapping is zero.
\end{lemma}
\begin{proof}
Take an element of $W_k[k-l_1,l_1+l_2]^{(0,N)}_{m,n,e}$ of the form
\bean
f_{i_1}\cdots f_{i_K}h_{j_1}\cdots h_{j_{K'}}v,
\eean
where $v$ is the highest weight vector of $W_k[k-l_1,l_1+l_2]$.
These vectors form a spanning set.

Suppose that $n-2m<l_1$, i.e., $K-K'<l_1$.
Let us prove that in $U_k[l_1,l_2]^{(0,N)}$ we have
\bean
f_{i_1}\cdots f_{i_K}h_{j_1}\cdots h_{j_{K'}}e_0^{l_1}u=0,
\eean
where the vector $u$ is the highest weight vector of $U_k[l_1,l_2]$.
Note that we have
\bean
e_0^{l_1+1}u=0.
\eean
This is the only property of $u$ we need for the proof.
We use induction on $K'$.
If $K'=0$ then $K<l_1$. Using
${\rm ad}\,e_0(f_i)=h_i$ and ${\rm ad}\,e_0(h_i)=0$, we have
\bean
({\rm ad}\,e_0)^{l_1}(f_{i_1}\cdots f_{i_K})u=0.
\eean
Therefore, we have
\bean
f_{i_1}\cdots f_{i_K}e_0^{l_1}u=
(-{\rm ad}\,e_0)^{l_1}(f_{i_1}\cdots f_{i_K})u=0.
\eean
If $K'>0$, we have
\bea
f_{i_1}\cdots f_{i_K}h_{j_1}\cdots h_{j_{K'}}e_0^{l_1}u
&=&
f_{i_1}\cdots f_{i_K}h_{j_1}\cdots h_{j_{K'-1}}e_0f_{j_{K'}}e_0^{l_1}u.
\label{EQ1}
\ena
Note that $K-(K'-1)<l_1+1$.
Changing $K'$ to $K'-1$, $l_1$ to $l_1+1$ and $u$ to $f_{j_{K'}}u$,
we apply the induction hypothesis.  Thus we obtain
\bea
f_{i_1}\cdots f_{i_K}h_{j_1}\cdots h_{j_{K'-1}}e_0^{l_1+1}f_{j_{K'}}u&=&0.
\label{EQ2}
\ena
{}From (\ref{EQ1}) and (\ref{EQ2}) we have
\bean
f_{i_1}\cdots f_{i_K}h_{j_1}\cdots h_{j_{K'}}e_0^{l_1}u
&=&
f_{i_1}\cdots f_{i_K}h_{j_1}\cdots h_{j_{K'-1}}e_0[f_{j_{K'}},e_0^{l_1}]u\\
&=&-l_1f_{i_1}\cdots f_{i_K}h_{j_1}\cdots h_{j_{K'}}e_0^{l_1}u\\
&=&0.
\eean
\end{proof}

\begin{lemma}\label{PROP2}
Suppose that $l_1+l_2\le k$.
The sequences
\bean
W_k^{(M,N)}[k-l_1,l_1+l_2]\rightarrow
U_k[l_1,l_2]^{(M,N)}\rightarrow
U_k[l_1-1,l_2]^{(M,N)}\rightarrow0\quad(M\ge0,N>0),\label{ONTO}\\
W_k^{(M,N)}[l_1+l_2,k-l_2]\rightarrow
V_k[l_1,l_2]^{(M,N)}\rightarrow
V_k[l_1,l_2-1]^{(M,N)}\rightarrow0\quad(M>0,N\ge0)\nonumber
\eean
are exact.
\end{lemma}
\begin{proof}
We consider the first sequence. Consider the following commutative diagram.
\[
\xymatrix{
&
&0\ar[d]
&0\ar[d]\\
&
&{\mathfrak{a}}^{(M,N)}U_k[l_1,l_2]\ar[r]\ar[d]
&{\mathfrak{a}}^{(M,N)}U_k[l_1-1,l_2]\ar[r]\ar[d]
&0\\
0\ar[r]
&\la e_0^{l_1}\ra\ar[r]\ar[d]
&U_k[l_1,l_2]\ar[r]^{\pi_1}\ar[d]^{\pi_2}
&U_k[l_1-1,l_2]\ar[r]\ar[d]
&0\\
0\ar[r]
&\pi_2\left(\la e_0^{l_1}\ra\right)\ar[r]\ar[d]
&U_k[l_1,l_2]^{(M,N)}\ar[r]\ar[d]
&U_k[l_1-1,l_2]^{(M,N)}\ar[r]\ar[d]
&0\\
&0&0&0}
\]
The columns and the first and second rows are exact.
One can check that the third row is exact by the standard diagram chasing.

Lemma \ref{PROP2} follows from the exactness of the third row 
and the surjectivity of \Ref{SH30}. 
It is clear from Lemma \ref{LEM} that the modification of the spaces
(\ref{MD2}) does not break the exactness.
\end{proof}

\subsection{Recursion of coinvariants}
In this section we establish a version of Proposition \ref{DECOM} for
the spaces of coinvariants which leads to 
a recursion for the space of coinvariants.

For $M,N\ge0$, consider the canonical surjection
\begin{equation}\label{PI}
\pi:W_k[l_1,l_2,l_3]\rightarrow W_k^{(M,N)}[l_1,l_2,l_3].
\end{equation}
We induce a filtration of $W_k^{(M,N)}[l_1,l_2,l_3]$ from (\ref{FILT2}):
\bea
W^{(M,N)}_{a,c}&=&\pi(W_{a,c}).
\label{MNFIL}
\ena
We call (\ref{MNFIL}) the canonical filtration of $W_k^{(M,N)}[l_1,l_2,l_3]$
and denote it by $C_{a,c}(W_k^{(M,N)}[l_1,l_2,l_3])$.

\begin{proposition}\label{COINVREC}
For the same $(a,c)$ and $(l_1',l_2')$ as in Proposition \ref{DECOM},
we have a surjection
\bea
T(W_k^{(M,N-1)}[l_1',l_2'])\rightarrow
{\rm Gr}^C_{a,c}(W_k^{(M,N)}[l_1,l_2,l_3]).
\label{MAP2}
\ena
\end{proposition}
\begin{proof}
First, we consider the spaces of coinvariants with respect to
${\mathfrak{a}}^{(M,N)}$ without the cutoff $(\ref{MD2})$.
We use the surjection
\bean
\bar\pi:W_k[l_1,l_2,l_3]\rightarrow
W_k[l_1,l_2,l_3]/{\mathfrak{a}}^{(M,N)}W_k[l_1,l_2,l_3].
\eean
Noting that $T({\mathfrak{a}}^{(M,N-1)})={\mathfrak{a}}^{(M,N)}$, we can induce
a surjection from (\ref{MAP}).
\bean
T(W_k[l_1',l_2']/{\mathfrak{a}}^{(M,N-1)}W_k[l_1',l_2'])
&\rightarrow&(W_{a,c}/W_{a,c+1})/{\mathfrak{a}}^{(M,N)}(W_{a,c}/W_{a,c+1})\\
&\simeq&W_{a,c}/(W_{a,c+1}+{\mathfrak{a}}^{(M,N)}W_{a,c}).
\eean
We continue to another surjection
\bean
&\rightarrow&W_{a,c}/
\left(W_{a,c+1}+{\mathfrak{a}}^{(M,N)}W_k[l_1,l_2,l_3]\cap W_{a,c}\right)\\
&\simeq&
(W_{a,c}+{\mathfrak{a}}^{(M,N)}W_k[l_1,l_2,l_3])/
(W_{a,c+1}+{\mathfrak{a}}^{(M,N)}W_k[l_1,l_2,l_3])\\
&\simeq&
\bar\pi(W_{a,c})/\bar\pi(W_{a,c+1}).
\eean
Therefore, we obtain the surjection
\bean
T(W_k[l_1',l_2']/{\mathfrak{a}}^{(M,N-1)}W_k[l_1',l_2'])
&\rightarrow&\bar\pi(W_{a,c})/\bar\pi(W_{a,c+1}).
\eean

Finally, we show that the cutoff (\ref{MD2}) does not break
the surjectivity.
There are two cases we must check: the case $(M,0)\rightarrow(M,1)$
and $(0,N-1)\rightarrow(0,N)$. In the first case, the weight spaces that
are cut off in (\ref{MD2}) are mapped to zero. 
The proof is similar to the proof of Lemma \ref{LEM}.

In the second case the proof is more elaborate. If $l_1+c\le k$,
the weight spaces that are cut off in the $(0,N-1)$-coinvariants,
are mapped to the weight spaces that are also cut off.
However, if $l_1+c>k$ (and therefore $l_1'=k-a$),
this is no longer true. In this case we show that
the weight spaces that are cut off are mapped to zero.

Let $K,K'\in\Z_{\ge0}$, and $\mu:\{1,\ldots,K+K'\}\rightarrow \Z$.
For a subset $J$ of $\{1,\ldots,K+K'\}$, we denote
\bean
f(\mu;J)=\prod_{j\in J}f_{\mu(j)},\quad h(\mu;J)=\prod_{j\in J}h_{\mu(j)}.
\eean
We set $w(\mu)=f(\mu;\{1,\ldots,K\})h(\mu;\{K+1,\ldots,K+K'\})$.

We claim that if $K-K'<a$, then 
\bea
w(\mu)h_1^af_1^cv=0\hbox{ in }W_{a,c}^{(0,N)}/W_{a,c+1}^{(0,N)}.
\label{EQW}
\ena 

First consider the case when $K<K'$. Since $e_0h_1^af_1^cv=0$ in
$W_{a,c}/W_{a,c+1}$, it is enough to show the equality
\bea
w(\mu)u\equiv0\bmod{\rm Im}\,e_0
\label{EQZERO}
\ena
in an arbitrary $\cH$-module assuming that $K'=K+1$ and $e_0u=0$.

Let $P_K$ be the set of all subsets of $\{1,\ldots,2K+1\}$
with the cardinality $K+1$, and denote by $J^c$ the complement of $J\in P_K$.
We have
\bean
[f(\mu;J),e_0]h(\mu;J^c)u\equiv0\bmod{\rm Im}\,e_0.
\eean
Define the matrix $M$ indexed by $P_K$:
\bean
M_{J,I}=
\begin{cases}
1\text{ if $I^c\subset J$};\\
0\text{ otherwise}.
\end{cases}
\eean
Note that $I^c\subset J$ is equivalent to $J=I^c\sqcup \{j\}$ for some $j\in J$
because $\sharp(I^c)=\sharp(J)-1$. Therefore,
$M$ is the incidence matrix of the term $f(\mu;I^c)h(\mu;I)u$
in $[f(\mu;J),e_0]h(\mu;J^c)u$.
It is easy to show that the matrix $M$ is non-degenerate
(e.g., by calculating its inverse). Therefore, (\ref{EQZERO}) follows.

If $K\ge K'$, we start with showing the equality
\bea
({\rm ad}\,e_0)^{K'}\bigl(f(\mu;\{1,\ldots,K+K'\})\bigr)h_1^af_1^cv\equiv0
\bmod{\rm Im}\,e_0.
\label{EQW2}
\ena
If $K'>0$, (\ref{EQW2}) follows from $e_0h_1^af_1^cv=0$ in $W_{a,c}/W_{a,c+1}$.
If $K'=0$, we have $K<a$. Since
$$
(c+1)Xh_1^af_1^cv=[X,e_0]h_1^{a-1}f_1^{c+1}v,
$$
(\ref{EQW2}) follows by induction on $K$.

Now we show (\ref{EQW}) by using (\ref{EQW2}).
For $\sigma\in{\mathfrak{S}}_{K+K'}$ we define
$\mu^\sigma:\{1,\ldots,K+K'\}\rightarrow \Z$ by $\mu^\sigma(n)=\mu(\sigma(n))$.
The left hand side of (\ref{EQW2}) is a positive linear combination of 
$w(\mu^\sigma)h_1^af_1^cv$
$(\sigma\in{\mathfrak{S}}_{K+K'})$. We will show that 
\bea
w(\mu^\sigma)h_1^af_1^cv=w(\mu)h_1^af_1^cv
\label{PERM}
\ena
for any $\sigma\in{\mathfrak{S}}_{K+K'}$. Then, our assertion (\ref{EQW})
follows from (\ref{EQW2}).

The symmetry (\ref{PERM}) is true for the transpositions
$1\leftrightarrow 2\leftrightarrow\cdots\leftrightarrow K$ and
$K+1\leftrightarrow K+2\leftrightarrow\cdots\leftrightarrow K+K'$.
Therefore, it is enough to show the symmetry with respect to
the transposition $1\leftrightarrow K+K'$.

We use the notations
\bean
(n)_m&=&(n+1)\cdots(n+m),\\
F^{(n)}&=&(-{\rm ad}\,e_0)^n\bigl(f(\mu;\{1,\ldots,K\})\bigr).
\eean
For $p,q\in\Z$ satisfying $0\leq p\leq K$, $0\leq q\leq K'-1$ and
$q\leq p\leq q+K-K'+1$,
we set
$$
[p,q;Q]=F^{(p)}\prod_{i\in Q}f_{\mu(i)}\prod_{i\in Q^c}h_{\mu(i)}
h_1^{a-p+q}f_1^{c+p-q}v.
$$
Here $Q$ is a subset $Q\subset\{K+1,\dots,K+K'\}$ such that $\sharp(Q)=q$,
and $Q^c$ is its complement in $\{K+1,\dots,K+K'\}$.
We have $[0,0;\emptyset]=w(\mu)h_1^af_1^cv$.
Note also that $p-q\leq K-K'+1\leq a$.

We will prove that $[p,q;Q]$ is in fact independent of $Q$ and
\begin{equation}
[p,q;Q]=(p-q)_q(c)_{p-q}[0,0;\emptyset].\label{EXACT}
\end{equation}
In particular, we have
$$
[K,K'-1;\{K+1,\dots,K+K'-1\}]=(K-K'+1)_{K'-1}(c)_{K-K'+1}[0,0;\emptyset],
$$
and therefore $w(\mu)h_1^af_1^cv$ is symmetric with respect to
the transposition $1\leftrightarrow K+K'$.

The proof is by induction on the pair of integers $(q,p)$
in the lexicographic ordering. For the induction steps we use
the identities
\begin{eqnarray}
&&Xh_1^af_1^cv=\frac1{(c)_r}(-{\rm ad}\,e_0)^r(X)
h_1^{a-r}f_1^{c+r}v\qquad(0\leq r\leq a),\label{I1}\\
&&F^{(p)}f_{j_1}\dots f_{j_q}h_{j_{q+1}}\dots h_{j_{K'}}h_1^{a-r}f_1^{c+r}v=
F^{(p+1)}f_{j_1}\dots f_{j_{q+1}}h_{j_{q+2}}\dots h_{j_{K'}}h_1^{a-r}f_1^{c+r}v
\nonumber\\
&&\quad-\sum_{1\leq s\leq q}F^{(p)}
\Bigl(\prod_{1\leq t\leq q\atop t\not=s}f_{j_t}\Bigr)
f_{j_{q+1}}h_{j_s}h_{j_{q+2}}\dots h_{j_{K'}}h_1^{a-r}f_1^{c+r}v
\nonumber\\
&&
+
\begin{cases}
0&\text{ if $r=0$};\\
-(c+r)F^{(p)}f_{j_1}\dots f_{j_{q+1}}h_{j_{q+2}}\dots h_{j_{K'}}
h_1^{a-r+1}f_1^{c+r-1}v&\text{ if $1\leq r\leq a$}.
\end{cases}
\label{I2}
\end{eqnarray}
Using (\ref{I1}) with $r=p$ we obtain
(\ref{EXACT}) with $q=0$. Because of the induction hypothesis, we can rewrite
(\ref{I2}) with $r=p-q$ as
$$
[p,q]=[p+1,q+1;Q]-q[p,q]
+
\begin{cases}
0&\text{ if $p=q$};\\
(c+p-q)[p,q+1]&\text{ otherwise}.
\end{cases}
$$
Here $[p,q]$ denotes the right hand side of (\ref{EXACT}) and $Q$ is an
arbitrary subset of $\{K+1,\dots,K+K'\}$ such that $\sharp(Q)=q+1$. From this
follows that $[p+1,q+1;Q]$ is independent of $Q$ and it is given by
(\ref{EXACT}).
\end{proof}

\begin{corollary}\label{RRR}
Suppose that for each $(a,c)$ satisfying
$0\leq a\leq l_3,0\leq c\leq l_2-a$ and $(l'_1,l'_2)$ given by (\ref{l primes})
we have a set of monomials $\{\m\}$ such that the set of vectors
$\{\m v_W[l'_1,l'_2,l'_1+l'_2-k]\}\subset W[l'_1,l'_2]$ forms a basis.
Then, the mapping
\bea\label{rrr}
\underset{\overset{\scriptstyle 0\le a\le l_3}{0\le c\le l_2-a}}
{\bigoplus}
W_k^{(M,N-1)}[l'_1,l'_2]
\rightarrow W_k^{(M,N)}[l_1,l_2,l_3]
\ena
which sends $\m v_W[l'_1,l'_2,l'_1+l'_2-k]$ to $T(\m)h_1^af_1^cv_W[l_1,l_2,l_3]$
is a surjection.
\end{corollary}

\subsection{Connection between Heisenberg and $\wsl$ modules}
\label{connection section}
In the course of the paper we will prove that the Heisenberg modules
${\rm Gr}^F(L_{k,l})$ and ${\rm Gr}^E(L_{k,l})$ are isomorphic to
$V_k[k-l,l]$ and $U_k[l,k-l]$, respectively. 
Now we are ready to establish surjectivity.
\begin{lemma}\label{PROS}
There exist surjective maps of $\cH$-modules,
\bea
V_k[k-l,l]&\rightarrow&{\rm Gr}^F(L_{k,l}),
\label{SURJ}\\
U_k[l,k-l]&\rightarrow&{\rm Gr}^E(L_{k,l}),
\label{SURJ2}
\ena
such that the highest weight vectors $v_V\in V_k[k-l,l]$
and $v_U\in U_k[l,k-l]$
are mapped to the image $v_F$ of $v_k[l]$ in ${\rm Gr}^F(L_{k,l})$ and to
the image $v_E$ of $u_k[l]$ in ${\rm Gr}^E(L_{k,l})$ respectively.
\end{lemma}
\begin{proof}
We consider the first map.
We need to show that the $\cH_k$-action on ${\rm Gr}^F(L_{k,l})$ satisfies
$e_1^{k-l+1}v_F=0$, $f_0^{l+1}v_F=0$ and $h_0v_F=0$. The first two
relations follow from the corresponding relations
for $v_k[l]\in L_{k,l}$. The last
one is proved in a similar manner as the proof of $[\bar h_i,\bar e_j]=0$
in Proposition \ref{IND}, by using the relation $h_0v_k[l]=lv_k[l]$
in $L_{k,l}$.
\end{proof}

\begin{corollary}\label{COINV COL}
For $M,N\geq 0$, there are surjective maps of coinvariants,
\bea
V_k[k-l,l]^{(M,N)}&\rightarrow&{\rm Gr}^F(L_{k,l})^{(M,N)},\label{FIN PROS}\\
U_k[l,k-l]^{(M,N)}&\rightarrow&{\rm Gr}^E(L_{k,l})^{(M,N)}.\label{FIN PROS2}
\ena
\end{corollary}
\begin{proof}
Corollary \ref{COINV COL} follows from Lemma \ref{PROS}.
\end{proof}

We will show that \Ref{SURJ}, \Ref{SURJ2} (and therefore \Ref{FIN PROS}, \Ref{FIN PROS2}) are 
isomorphism in Corollary \ref{GRA}.

In order to construct a recursion we need to "decompose" the spaces of the
coinvariants into smaller pieces.

\begin{lemma}\label{SPLIT}
There exist surjective linear maps
\bean
\bigoplus\limits_{i=0}^lW_k[k-l+i,k-i]&\to& V_k[k-l,l],\\
\bigoplus\limits_{i=0}^lW_k[k-l+i,k-i]&\to& U_k[l,k-l].\\
\eean
Moreover, for $M,N\geq 0$, 
these maps can be factorized to the surjective maps of coinvariants,
\begin{eqnarray}
\bigoplus\limits_{i=0}^lW^{(M,N)}_k[k-l+i,k-i]
&\to& V_k[k-l,l]^{(M,N)} \qquad (M>0),\label{VMN}\\
\bigoplus\limits_{i=0}^lW^{(M,N)}_k[k-l+i,k-i]
&\to& U_k[l,k-l]^{(M,N)} \qquad (N>0).\label{UMN}
\end{eqnarray}
\end{lemma}
\begin{proof}
Lemma \ref{SPLIT} follows from Corollary \ref{TRIPLES} and Lemma \ref{PROP2}.
\end{proof}

We will see below (Corollary \ref{SES cor}) that all these maps are, in fact,
isomorphisms.

We note that if $0<l\leq k$ and $M=0$ the dimension of
$V[k-l,l]^{(0,N)}$ is greater than the dimension of
$L_{k,l}^{(0,N)}$. This can be seen from the following lemma. 

\begin{lemma}\label{ZERO}
For $s>-l$,
\bean
\left(L^{(0,N)}_{k,l}\right)_{s,e}=0.
\eean
\end{lemma}
\begin{proof}
The $\widehat{\mathfrak{sl}}_2$-module $L_{k,l}$ is generated by
the vector $u_k[l]=f_0^lv_k[l]$. The set of vectors of the form
\bean
f_{i_1}\cdots f_{i_K}h_{j_1}\cdots h_{j_{K'}}u_k[l]\in
\left(L^{(0,N)}_{k,l}\right)_{-2K-l,e}
\eean
is a spanning set of $L_{k,l}^{(M,N)}$. In fact, the operators $e_i$
can be eliminated either by using the highest weight condition
$e_iu_k[l]=0$ $(i\le-1)$ or by taking the quotient with respect to
the subspace ${\rm Im}\,e_i$ $(i\ge0)$. The statement of the lemma follows.
\end{proof}

The weight space $V_k[k-l,l]^{(0,N)}_{m,n,e}$ is mapped to
${\rm Gr}^F(L_{k,l})^{(M,N)}_{l+2m-2n,e}$ by (\ref{SURJ}). If $l+2m-2n>-l$,
by Lemma \ref{ZERO} the space ${\rm Gr}^F(L_{k,l})^{(0,N)}_{l+2m-2n,e}$ 
is zero. It happens, in particular, for $l>0$ and
$m=n=0$. On the other hand, the space
$V_k[k-l,l]^{(0,N)}_{0,0}$ is spanned by the the highest weight vector and is
one dimensional.

Note that in the case $N>0$, even if we have $M=0$, we can use the
$U_k[l,k-l]^{(M,N)}$ instead. If $M=N=0$, the space of
$\wsl$-coinvariants is easily calculated. The result is
\bean
L^{(0,0)}_{k,l}&=&
\begin{cases}
\C \cdot v_k[0]\text{ if $l=0$;}\\
0\text{ otherwise}.
\end{cases}
\eean

\subsection{Dimension counting}\label{zipper section}
In this section we show that the recursion for the spaces
$W^{(M,N)}[l_1,l_2]$ coincides with the recursion for the combinatorial paths
${\mathcal C}_{k,l}^{(M+N)}$.
This gives us the upper bound of ${\rm dim}\,L^{(M,N)}_{k,l}$
by $d^{(M+N)}_{k,l}$. Since we have already
shown that the same number gives the lower bound, we complete
the proof of the equality ${\rm dim}\,L^{(M,N)}_{k,l}=d^{(M+N)}_{k,l}$. 

We start with
\begin{lemma}\label{PROP1}
\bean
{\rm dim}\,V_k[k-l,l]^{(M,N)}\geq {\rm dim}\,L^{(M,N)}_{k,l}\geq d_{k,l}^{(M+N)},\\
{\rm dim}\,U_k[l,k-l]^{(M,N)}\geq {\rm dim}\,L^{(M,N)}_{k,l}\geq d_{k,l}^{(M+N)}.
\eean
\end{lemma}
\begin{proof}
We show the first line of inequalities, the second one is proved similarly.

Using Lemma \ref{deformation} we have
\bean
{\rm dim}\,{\rm Gr}^F(L_{k,l})^{(M,N)}\geq {\rm dim}\,L^{(M,N)}_{k,l}.
\eean
{}From Lemma \ref{PROS} we have
\bean
{\rm dim}\,V_k[k-l,l]^{(M,N)}\geq {\rm dim}\,{\rm Gr}^F(L_{k,l})^{(M,N)}.
\eean
This lemma follows from these inequalities along with
Theorem \ref{Ver} and Proposition \ref{upper bound}.
\end{proof}

\begin{corollary}\label{dim ineq cor}
\be
\sum_{i=0}^l\dim (W^{(M,N)}[k-l+i,k-i])\geq d_{k,l}^{(M+N)}.
\ee
\end{corollary}
\begin{proof}
This corollary follows from Lemmas \ref{PROP1} and \ref{SPLIT}. 
\end{proof}

Now we show that the recursive relation described by Proposition
\ref{COINVREC} for the spaces of coinvariants $W^{(M,N)}[l_1,l_2]$
coincides with the recursion of the combinatorial sets
${\mathcal C}_{k,l}^{(N)}$
(see \Ref{combinatorial recursion} and Proposition \ref{comb rec theorem}). 

\begin{lemma}\label{LEMREC}
The set of
${\mathcal C}_{k,l'}^{(N)}[i']$ which appear in the summation of the right hand side of
(\ref{comb recur}) coincides with the set of $W_k^{(M,N-1)}[l'_1,l'_2]$
which appear in the summation of the left hand side of (\ref{rrr}) in the case
$l_1=k-l+i,l_2=k-i,l_3=k-l$ by the identification of
${\mathcal C}_{k,l'}^{(N)}[i']$ with  $W_k^{(M,N-1)}[k-l'+i',k-i']$.
\end{lemma}
\begin{proof}
Since $l''=2i-l+l'$, $l'+l''\geq l$ and $l+l'\geq l''$ are equivalent to
$i\geq0$ and $i\leq l$, respectively.
Therefore, the admissibility of the triple $(l,l'',l')$ reduces to
\be
l+l'+l''\leq 2k, \qquad l''+l'-l\geq 0.
\ee
Namely, the sum in \Ref{comb recur} is taken over $l'$ and $i'$ such that 
\be
l-i\leq l'\leq k-i,\qquad 0\leq i' \leq l'.
\ee
This is equal to the sum over
\begin{equation}
\label{SUM1}
0\leq i'\leq k-i,\qquad \max(l-i,i')\leq l'\leq k-i.
\end{equation}

Now consider the surjection
\bea\label{SURW}
\underset{\overset{\scriptstyle 0\le a\le k-l}{0\le c\le k-i-a}}
{\bigoplus}
W_k^{(M,N-1)}[l'_1,l'_2]
\rightarrow W_k^{(M,N)}[k-l+i,k-i]
\ena
The sum in the left hand side is equivalent to the sum over
\begin{equation}
\label{SUM2}
0\leq c\leq k-i,\qquad 0\leq a \leq \min(k-l,k-i-c).  
\end{equation}
Set
\begin{eqnarray}
i'&=&c,\label{IDENT-a}\\
l'&=&
\begin{cases}
a+l-i\quad&\text{if $0\leq c\leq l-i$};\\
a+c\quad&\text{if $l-i\leq c\leq k-i$}.
\end{cases}
\label{IDENT-b}
\end{eqnarray}
Then, we have
\begin{equation}
l'_1=k-l'+i',\quad l'_2=k-i'.\label{IDENT2}
\end{equation}
The sum (\ref{SUM1}) is equal to (\ref{SUM2}) by this identification.
\end{proof}

Now we are ready to prove
\begin{theorem}\label{MAINTH}
For $M,N\geq 0$, 
\begin{equation}\label{mainth}
\dim L_{k,l}^{(M,N)}=d^{(M+N)}_{k,l}.
\end{equation}
\end{theorem}
\begin{proof}
We will prove the above equality by induction on $(M,N)$ in three kinds
steps. The first steps are $(0,N)\rightarrow(0,N+1)$, the second
$(M,0)\rightarrow(M,1)$, and the third $(M,N)\rightarrow(M,N+1)$.

In the first and the third steps, we assume that
\be 
\dim W_k^{(M,N-1)}[k-l+i,k-i]=\sharp({\mathcal C}_{k,l}^{(M+N-1)}[i]),
\ee
and prove that
\begin{equation}
\dim W_k^{(M,N)}[k-l+i,k-i]=\sharp({\mathcal C}_{k,l}^{(M+N)}[i]).
\label{AUX}
\end{equation}
In the second steps, we assume that
\begin{equation}
\dim W_k^{(M,0)}[k-i,k-l+i]=\sharp({\mathcal C}_{k,l}^{(M)}[i]),
\label{2ASS}
\end{equation}
and prove that
\be 
\dim W_k^{(M,1)}[k-l+i,k-i]=\sharp({\mathcal C}_{k,l}^{(M+1)}[i]).
\ee
In each step, we also prove the equality of the form (\ref{mainth}).
The base of the induction is $(M,N)=(0,0)$, where (\ref{AUX}) is valid
because of (\ref{Z}) and (\ref{INIT}). Note also that the assumption
(\ref{2ASS}) follows from (\ref{AUX}) with $(M,N)$ replaced with $(0,M)$,
which is proved in the first steps.

Now, we show the first and the third induction steps at the same time.

Using Theorem \ref{Ver}, Proposition \ref{upper bound},
Corollary \ref{COINV COL} and Lemma \ref{SPLIT},
we have the chain of inequalities:
$$
d^{(M+N)}_{k,l}
=\dim L_{k,l}^{(M,N)}({\bf z},{\bf z'})
\leq \dim L_{k,l}^{(M,N)}
\leq \sum _{i=0}^l \dim W_k^{(M,N)}[k-l+i,k-i].
$$

Using (\ref{SURW}) we obtain
\begin{eqnarray}
\dim W^{(M,N)}_k[k-l+i,k-i]&\leq&\nonumber\\
\sum_{l-i\leq l'\leq k-i\atop0\leq i' \leq l'}
\dim W^{(M,N-1)}_k[k-l'+i',k-i']
&=&\sum_{l-i\leq l'\leq k-i\atop0\leq i' \leq l'}
\sharp({\mathcal C}^{(M+N-1)}_{k,l'}[i'])
\label{DIF}\\
&=&\sharp({\mathcal C}_{k,l}^{(M+N)}[i]).\nonumber
\end{eqnarray}

Summing up these inequalities for $i=0,\dots,l$ and finally using
Corollary \ref{CARLEM} we obtain 
\be
\sum_{i=0}^l\dim W^{(M,N)}_k[k-l+i,k-i]\leq
\sum_{i=0}^l \sharp ({\mathcal C}_{k,l}^{(M+N)}[i])
=\sharp({\mathcal C}_{k,l}^{(M+N)})=d^{(M+N)}_{k,l}.  
\ee
In particular, we obtain (\ref{AUX}) and (\ref{mainth}).

In the second steps we proceed with $N=1$.
We only append one formula to (\ref{DIF}):
$$
\sum_{l-i\leq l'\leq k-i\atop0\leq i' \leq l'}
\dim W^{(M,0)}_k[k-l'+i',k-i']
=\sum_{l-i\leq l'\leq k-i\atop0\leq i' \leq l'}
\dim W^{(M,0)}_k[k-i',k-l'+i']
=\sum_{l-i\leq l'\leq k-i\atop0\leq i' \leq l'}
\sharp({\mathcal C}^{(M)}_{k,l'}[i'])
$$
The rest of proof is similar.
\end{proof}

\begin{corollary}\label{LEMSUR}
The maps (\ref{MAP}) and (\ref{MAP2}) are isomorphisms if $l_3=l_1+l_2-k$.
\end{corollary}
\begin{proof}
The corollary follows from the proof of Theorem \ref{MAINTH}.
\end{proof}

In Section 6, we will prove that (\ref{rrr}) is an isomorphism
in all other cases, too.

\begin{corollary}\label{CORL}
Suppose $N>0$. We have the equality
\begin{equation}\label{NOJUMP}
{\rm dim}\,L^{(M,N)}_{k,l}={\rm dim}\,{\rm Gr}^E(L_{k,l})^{(M,N)}.
\end{equation}
The mapping
\begin{equation}
\bigoplus\limits_{i=0}^lW^{(M,N)}_k[k-l+i,k-i]
\to{\rm Gr}^E(L_{k,l})^{(M,N)}\label{INV}
\end{equation}
obtained by the composition of (\ref{UMN}) and (\ref{FIN PROS2}),
is an isomorphism.
\end{corollary}
\begin{proof}
This is proved in the proof of Theorem \ref{MAINTH}.
\end{proof}

\begin{corollary} \label{GRA}
The maps \Ref{SURJ}, \Ref{SURJ2} are isomorphisms of $\cH$ modules.
\end{corollary}
\begin{proof}
The corollary follows from the proof of Theorem \ref{MAINTH}.
\end{proof}

\begin{corollary}\label{SES cor}
For $l_1,l_2$ such that $0\le l_1,l_2\le k\le l_1+l_2$,
we have the short exact sequences
\bean
0\rightarrow W_k^{(M,N)}[l_1+l_2,k-l_2]\rightarrow V_k[l_1,l_2]^{(M,N)}
\rightarrow V_k[l_1,l_2-1]^{(M,N)}\rightarrow0 \qquad (M>0,N\geq 0),\\
0\rightarrow W_k^{(M,N)}[k-l_1,l_1+l_2]\rightarrow U_k[l_1,l_2]^{(M,N)}
\rightarrow U_k[l_1-1,l_2]^{(M,N)}\rightarrow0 \qquad (N\geq 0,M> 0).
\eean
\end{corollary}
\begin{proof}
By Lemma \ref{PROP2}, all we have to show is the injectivity of
the mappings
\bean
W_k^{(M,N)}[l_1+l_2,k-l_2]\rightarrow V_k[l_1,l_2]^{(M,N)},\\
W_k^{(M,N)}[k-l_1,l_1+l_2]\rightarrow U_k[l_1,l_2]^{(M,N)}.
\eean
The injectivity follows from the proof of Theorem \ref{MAINTH}.
\end{proof}

\begin{corollary}\label{SES cor2}
For $l_1,l_2$ such that $0\le l_1,l_2\le k\le l_1+l_2$,
we have the short exact sequences
\bean
0\rightarrow W_k[l_1+l_2,k-l_2]\rightarrow V_k[l_1,l_2]
\rightarrow V_k[l_1,l_2-1]\rightarrow0,\\
0\rightarrow W_k[k-l_1,l_1+l_2]\rightarrow U_k[l_1,l_2]
\rightarrow U_k[l_1-1,l_2]\rightarrow0.
\eean
\end{corollary}

\begin{proof}
This follows from Corollary \ref{SES cor} by letting $M,N\rightarrow\infty$.
\end{proof}

\begin{corollary}\label{FHBASIS}
The set of monomial vectors
\bea\label{FHBAS}
\{f_{N-1}^{a_{N-1}} h_{N-1}^{b_{N-1}}\dots
f_1^{a_1}h_1^{b_1}f_0^{a_0}v_k[l]\in L_{k,l}^{(0,N)}({\bf z})
\; ; \;({\bf a};{\bf b})
\in{\mathcal C}^{(N)}_{k,l}\}
\ena
forms a basis of $L_{k,l}^{(0,N)}({\bf z})$ for all ${\bf z}\in \C^N$.
\end{corollary}
\begin{proof}
The recursion (\ref{rrr}) gives us a recursive way of constructing
monomial basis of the space of coinvariants $W_k^{(M,N)}[l_1,l_2]$.
As shown in Lemma \ref{LEMREC}, for each $(i,l)$ satisfying $0\leq i\leq l$,
the mapping (\ref{IDENT-a}), (\ref{IDENT-b}) gives a bijection between the set
of $(a,c)$ satisfying $0\leq a\leq k-l,0\leq c\leq k-i-a$ and the set of
$(i',l')$ satisfying $0\leq i'\leq l'$ and such that $(l,l''=2i-l+l',l')$
is admissible. Recall also that under this identification,
$(l'_1,l'_2)$ given by (\ref{IDENT2}) satisfies (\ref{l primes})
with $l_1=k-l+i$.

Suppose that we have a set of monomials $\{\m\}$ for each $(i',l')$ such that
for each $l'$ the union $\sqcup_{0\leq i'\leq l'}\{\m f_0^{i'}v_k[l']\}$ 
forms a basis of
${\rm Gr}^E(L_{k,l'})^{(M,N-1)}$. The inverse map of the bijection
(\ref{INV}) (with $N,l,i$ replaced by $N-1,l',i'$) maps the basis
$\{\m f_0^{i'}v_k[l']\}$ to the basis
$\{\m v_W[l'_1,l'_2,l'_1+l'_2-k]\}\subset W_k^{(M,N-1)}[l'_1,l'_2]$.
By Corollary \ref{RRR} we have the basis $\{T(\m)h_1^af_1^cv_w[k-l+i,k-i,k-l]\}$
of $W_k^{(M,N)}[k-l+i,k-i]$. By (\ref{INV}) we then have the basis
$\{T(\m f_0^{i'})h_1^af_1^cv_k[l]\}$ of ${\rm Gr}^E(L_{k,l})^{(M,N)}$.

In this way we obtain the basis of ${\rm Gr}^E(L_{k,l})^{(M,N)}$ from
the basis of ${\rm Gr}^E(L_{k,l'})^{(M,N-1)}$.
Noting that $i=(l+l''-l')/2,a=i+l'-l-(i'-l+i)^+$
and comparing these formulas with (\ref{RCOM1}), (\ref{RCOM2})
we see that this recursion gives exactly the basis (\ref{FHBAS}).
\end{proof}

\begin{example}
As an example, we list
the $k=1$ monomials for $M+N\leq3$:
\bean
\begin{array}{cccc}
(M,N)& W_1^{(M,N)}[1,1]& W_1^{(M,N)}[1,0]& W_1^{(M,N)}[0,1]\\
(0,0)&1&\emptyset&\emptyset\\
(0,1)&1&1&\emptyset\\
(1,0)&1&\emptyset&1\\
(0,2)&1,f_1&1&f_1\\
(1,1)&1,h_1&1&1\\
(2,0)&1,e_1&e_1&1\\
(0,3)&1,f_1,f_2,f_2h_1&1,f_2&f_1,f_2\\
(1,2)&1,f_1,h_1,h_2&1,h_2&1,f_1\\
(2,1)&1,e_1,h_1,h_2&1,e_1&1,h_2\\
(3,0)&1,e_1,e_2,e_2h_1&e_1,e_2&1,e_2
\end{array}
\eean
\end{example}

\begin{remark}
For small values of $(M,N)$ we have the following results.
(Here, we assume that $l_3\leq{\rm min}(l_1,l_2)$.)
\begin{eqnarray*}
W^{(0,0)}[l_1,l_2,l_3]&=&
\begin{cases}
\C\cdot1&\text{ if $l_1=l_2=k$;}\\
0&\text{ otherwise,}
\end{cases}\\
W^{(0,1)}[l_1,l_2,l_3]&=&
\begin{cases}
\C\cdot1&\text{ if $l_1=k$;}\\
0&\text{ otherwise,}
\end{cases}\\
W^{(1,0)}[l_1,l_2,l_3]&=&
\begin{cases}
\C\cdot1&\text{ if $l_2=k$;}\\
0&\text{ otherwise,}
\end{cases}\\
W^{(1,1)}[l_1,l_2,l_3]&=&
\oplus_{0\leq a \leq l_3}\C h_1^a.
\end{eqnarray*}
Note also that if we are interested in the case $M,N>0$,
we can avoid the argument on cut-off and apply the recursion relation
(\ref{rrr}) with the base $(M,N)=(1,1)$.
\end{remark}
\section{Miscellaneous results}\label{SEC6}
In this section we give several results which follow from the recursion of
coinvariants (\ref{rrr}). Some of them will play an important role in our next
paper where we will compute the characters of the coinvariants.
\subsection{Generalized recursion}
\label{general recur section}
In this section we prove the decomposition of the space of coinvariants
$W^{(M,N)}_k[l_1,l_2,l_3]$ for generic $l_3$. For this purpose we need
to introduce another filtration of $W_k[l_1,l_2,l_3]$. We will show that the
adjoint graded spaces of the two filtrations are canonically isomorphic.

The filtration we used in Section \ref{SEC-COIN}, i.e.,
the canonical filtration of the first kind, is such that we take the first
filtration with respect to $h_1^a$ and then the second one with respect to
$f_1^c$. We will define the filtration in the opposite
order. Namely, we consider the subspaces
$$
C'_{a,c}(W_k[l_1,l_2,l_3])=\la h_1^af_1^c,f_1^{a+c+1}\ra\subset
W_k[l_1,l_2,l_3]
$$
for $(a,c)$ satisfying (\ref{REGION}). This region is equivalently written as
$$
0\leq c\leq l_2,\quad 0\leq a\leq{\rm min}(l_3,l_2-c).
$$
We abbreviate $C'_{a,c}(W_k[l_1,l_2,l_3])$ to $W'_{a,c}$.
We define $l'_1,l'_2$ as before, i.e., by (\ref{l primes}).

Note the inclusion $\la h_1^af_1^{b-a}\ra\supset\la h_1^{a+1}f_1^{b-a-1}\ra$
is valid. We have the filtration:
\begin{eqnarray*}
&&W_k[l_1,l_2,l_3]=W'_{0,0}\supset\dots\\
&&\quad\supset W'_{0,c}\supset W'_{1,c-1}\supset
W'_{{\rm min}(l_3,c),c-{\rm min}(l_3,c)}
\supset W'_{0,c+1}\\
&&\qquad\qquad\qquad\qquad\qquad\qquad\qquad\qquad\qquad\qquad\qquad
\dots\supset W'_{l_3,l_2-l_3}\supset W'_{0,l_2+1}=0.
\end{eqnarray*}

The arguments in Sections \ref{SEC3} and \ref{SEC-COIN}
can be repeated for the canonical filtration of the second kind.
We describe the results without giving proofs. The proofs can be
repeated without modifications by the following reason.

Set $u_0=h_1^af_1^cv_W\in W_k[l_1,l_2,l_3]$,
$u_1=h_1^af_1^cv_W\in\la h_1^af_1^c,h_1^{a+1}\ra
/\la h_1^af_1^{c+1},h_1^{a+1}\ra$
and $u_2=h_1^af_1^cv_W\in\la h_1^af_1^c,f_1^{a+c+1}\ra
/\la h_1^{a+1}f_1^{c-1},f_1^{a+c+1}\ra$.
In the proofs, we use additional properties of $u_1$ or $u_2$,
which are not satisfied by $u_0$. Namely, we use
$$
e_0u_i=0,\quad f_1u_i=0\quad(i=1,2).
$$
These are equivalent to $h_1^af_1^{c+1}v_W=0$ and $h_1^{a+1}f_1^{c-1}v_W=0$,
respectively. They are valid in both of the quotient spaces.

\begin{proposition}\label{DECOM-2}
For each $(a,c)$ satisfying $(\ref{REGION})$ there exists an $\cH$-linear
surjection
\begin{eqnarray}
T(W_k[l_1',l_2'])&\rightarrow&{\rm Gr}^{C'}_{a,c}(W_k[l_1,l_2,l_3]),
\label{MAP-2}\\
T(mv_W[l'_1,l'_2,l'_1+l'_2-k])&\mapsto&T(m)v_{a,c},\nonumber\\
v_{a,c}&=&h_1^af_1^cv_W[l_1,l_2,l_3].\nonumber
\end{eqnarray}
Here $l'_1,l'_2$ are given by $(\ref{l primes})$.
\end{proposition}

\begin{proposition}\label{COINVREC-2}
For the same $(a,c)$ and $(l_1',l_2')$ as in Proposition \ref{DECOM-2},
we have a surjection
\bea
T(W_k^{(M,N-1)}[l_1',l_2'])\rightarrow
{\rm Gr}^{C'}_{a,c}(W_k^{(M,N)}[l_1,l_2,l_3]).
\label{MAP2-2}
\ena
\end{proposition}

\begin{proposition}\label{PRE}
If $l_3=l_1+l_2-k$, the maps $(\ref{MAP-2})$ and $(\ref{MAP2-2})$
are isomorphisms.
\end{proposition}

Now, we remove the above restriction on $l_3$.
\begin{proposition}\label{RECURSION}
The maps $(\ref{MAP-2})$ and $(\ref{MAP2-2})$ are isomorphisms.
\end{proposition}
\begin{proof}
Consider the exact sequence obtained by the same argument as in
the proof of Lemma \ref{PROP2}.
$$
0\rightarrow\pi(\la h_1^{l_3+1}\ra)\rightarrow W^{(M,N)}_k[l_1,k]
\rightarrow W^{(M,N)}_k[l_1,k,l_3]\rightarrow0.
$$
Here $\pi:W_k[l_1,k]\rightarrow W^{(M,N)}_k[l_1,k]$ is the canonical surjection
and $\la h_1^{l_3+1}\ra$ is the submodule of $W_k[l_1,k]$
generated by $h_1^{l_3+1}v_W[l_1,k,l_1]$.

The subspace $\la h_1^{l_3+1}\ra$ appears in the canonical filtration of the
first kind as $W_{l_3+1,0}$. Therefore, by Corollary \ref{LEMSUR}. we know
the exact decomposition of the adjoint graded spaces for 
$\pi(\la h_1^{l_3+1}\ra)$ and $W_k^{(M,N)}[l_1,k]$. In particular, we have
$$
{\rm dim}\,\pi(\la h_1^{l_3+1}\ra)=
\sum_{l_3+1\leq a\leq l_1\atop0\leq c\leq k-a}
{\rm dim}\,W_k^{(M,N-1)}[l'_1,l'_2]
$$
and
$$
{\rm dim}\,W_k^{(M,N)}[l_1,k]=
\sum_{0\leq a\leq l_1\atop0\leq c\leq k-a}
{\rm dim}\,W_k^{(M,N-1)}[l'_1,l'_2].
$$
Therefore, we have
$$
{\rm dim}\,W_k^{(M,N)}[l_1,k,l_3]=
\sum_{0\leq a\leq l_3\atop0\leq c\leq k-a}
{\rm dim}\,W_k^{(M,N-1)}[l'_1,l'_2].
$$
{}From this the statement of the proposition for $l_2=k$ follows.

Next, consider the exact sequence
$$
0\rightarrow\pi(\la f_1^{l_2+1})\ra\rightarrow W^{(M,N)}_k[l_1,k,l_3]
\rightarrow W^{(M,N)}_k[l_1,l_2,l_3]\rightarrow0.
$$
The subspace $\la f_1^{l_2+1}\ra$
appears in the canonical filtration of the second kind
as $W'_{0,l_2+1}$. Applying the same argument as before, we obtain
the statement of the proposition for general $l_2$.
\end{proof}

Let us summarize the results obtained in
Propositions \ref{DECOM}, \ref{COINVREC} and \ref{RECURSION}.

\begin{theorem}\label{R}
Let $l_1,l_2,l_3$ be such that
$0\leq l_1,l_2\leq k$, $0\leq l_3\leq\min(l_1,l_2)$.
There exists an isomorphism of vector spaces
\bean
\rho:
\bigoplus\limits_{a,c} W_k[l_1',l_2']
\rightarrow W_k[l_1,l_2,l_3]
\eean
where 
\bean
l_1'=l_1+c-a-(l_1+c-k)^+,\quad l_2'=k-c,
\eean 
and the sum in the left hand side is taken over
$0\leq a\leq l_3,\quad0\leq c\leq l_2-a.$
This map is explicitly described in Proposition \ref{DECOM}. Moreover, the induced map
\bean
\rho_k^{(M,N)}:\bigoplus\limits_{a,c} W_k^{(M,N-1)}[l_1',l_2']
\rightarrow W_k^{(M,N)}[l_1,l_2,l_3]
\eean
is also an isomorphism of vector spaces.
\end{theorem}

It is convenient to write the above theorem informally:
\bean
W_k^{(M,N)}[l_1,l_2,l_3]&=&\bigoplus\limits_{\stackrel{0\leq a\leq
l_3}{0\leq c\leq l_2-a}}T(W_k^{(M,N-1)}[l_1',l_2'])h_1^af_1^c.
\eean

Note that because the roles of $e_i$ and $f_i$ in $W$ modules are similar, $\iota (W_k[l_1,l_2])=W_k[l_2,l_1]$, 
we also have a recursion relations on $M$. Namely, 
\bean
W_k^{(M,N)}[l_2,l_1,l_3]&=&\bigoplus\limits_{\stackrel{0\leq a\leq
l_3}{0\leq c\leq l_2-a}}T'(W_k^{(M-1,N)}[l_2',l_1'])h_1^ae_1^c,
\eean
where $T'$ is such the automorphism of $\cH$ that $\iota T=T'\iota$, namely
\begin{equation}\label{T'}
T'(e_i)=e_{i+1},\qquad T'(f_i)=f_i,\qquad T(h_i)=h_{i+1}.
\end{equation}

\begin{remark} In fact we have constructed many bases of the $(M_0,N_0)$
coinvariants. Namely, we start from $(0,0)$ case and choose a path on the
$(M,N)$ plane from $(0,0)$ to $(M_0,N_0)$, increasing either $M$ or $N$ by one.
Then we apply the recursions in the order according to this path and get a
basis. Clearly, the bases corresponding to different paths are different.
\end{remark}

\subsection{Short exact sequences}
In this section we establish several short exact sequences. 

We have the following exact sequences:
\bea
&&0\rightarrow \la h_0\ra\rightarrow \bV_k[l_1,l_2]
\rightarrow V_k[l_1,l_2]\rightarrow0\quad(l_1+l_2\leq k),\label{52}\\
&&0\rightarrow \la h_0\ra\rightarrow \bU_k[l_1,l_2]
\rightarrow U_k[l_1,l_2]\rightarrow0\quad(l_1+l_2\leq k),\\
&&0\rightarrow \la h_1^{l_1+l_2-k+1}\ra\rightarrow
\bW_k[l_1,l_2]\rightarrow W_k[l_1,l_2]\rightarrow0,\quad(l_1+l_2\geq k).
\label{53}
\ena

\begin{proposition}\label{triple prop2}
The following are well-defined and surjective $\cH$-linear mappings.
\bea
\bV_k[l_1-1,l_2-1]&\rightarrow& \la h_0\ra\subset \bV_k[l_1,l_2]
\quad(l_1+l_2\leq k),
\label{SH1}\\
\bU_k[l_1-1,l_2-1]&\rightarrow& \la h_0\ra\subset \bU_k[l_1,l_2]
\quad(l_1+l_2\leq k),
\label{SH10}\\
\bW_k[k-l_2-1,k-l_1-1]&\rightarrow& \la h_1^{l_1+l_2-k+1}\ra
\subset \bW_k[l_1,l_2]
\quad(l_1+l_2\geq k).\label{SH2}
\ena
Each of these mappings is defined to be such that the highest weight vector of
the left hand side is mapped to the indicated generator of the submodule in the right hand side.
\end{proposition}
\begin{proof}
We have to show the validity of the relations (\ref{HIGH})
and (\ref{REL1}), or (\ref{REL2}) and (\ref{REL3}).
We will give the proofs for the non-trivial ones.
The relations $e_1^{l_1}h_0v=f_0^{l_2}h_0v=0$ for
(\ref{SH1}) follow from Lemma \ref{LM0} with $c=0$.
The relations $e_1^{k-l_2}h_1^{l_1+l_2-k+1}v=f_1^{k-l_1}h_1^{l_1+l_2-k+1}v=0$
for (\ref{SH2}) follow from Lemma \ref{LM2} with $c=0$.
\end{proof}

\begin{proposition}\label{tr prop}
For $l_1,l_2$ such that $0\le l_1\leq l_2\le k$ and $l_1+l_2\ge k$, we have the
exact sequence
\bean
0\rightarrow W_k^{(M,N)}[k-l_2-1,k-l_1-1,k-l_2-1]\rightarrow
W_k^{(M,N)}[l_1,l_2,l_1]\rightarrow W_k^{(M,N)}[l_1,l_2]
\rightarrow0.
\eean
\end{proposition}
\begin{proof}
The proof goes similarly as in Proposition \ref{RECURSION} except for the
difference in the formulas of $l'_1,l'_2$ between the first term and
the second and third terms.

The difference of the ranges of $(a,c)$ for
$W_k^{(M,N)}[l_1,l_2,l_1]$ and $W_k^{(M,N)}[l_1,l_2,l_1+l_2-k]$ is given by
\bean
l_1+l_2-k+1\leq a\leq l_1,0\leq c\leq l_2-a.
\eean
In this range, we have
\bean
l_1+c-k<0.
\eean
Therefore, we have
\bean
l_1'=l_1+c-a,l_2'=k-c.
\eean
On the other hand, the range of $(\bar a,\bar c)$ for
$W_k^{(M,N)}[k-l_2-1,k-l_1-1,k-l_2-1]$ is
\bean
0\leq\bar a\leq k-l_2-1,0\leq \bar c\leq k-l_1-1-\bar a.
\eean
Noting that $l_1+l_2\geq k$ we have
\bean
k-l_2-1+\bar c-k<0.
\eean
Therefore, we have
\bean
l_1'=k-l_2-1+\bar c-\bar a,l_2'=k-\bar c.
\eean
By the correspondence
\bean
a=\bar a+l_1+l_2-k+1,c=\bar c,
\eean
these two ranges of $(l_1',l_2')$ coincide.
\end{proof}

\begin{lemma}
For $l_1,l_2$ such that $0\le l_1,l_2\le k$ and $l_1+l_2\le k$ we have an exact sequence
\bea\label{triple6}
0\rightarrow W_k^{(M,N)}[l_1-1,l_2-1,l_3-1]\rightarrow W_k^{(M,N)}[l_1,l_2,l_3]
\rightarrow W_k^{(M,N)}[l_1,l_2,0]\rightarrow0,
\ena
where $l_3=\min(l_1,l_2)$.
\end{lemma}
\begin{proof}
The proof goes similarly as in the previous proposition.

Noting that $l_1+l_2\leq k$, we have
\bea
l_1'=l_1+c-a,l_2'=k-c,
\label{RANGE}
\ena
for $W_k^{(M,N)}[l_1,l_2,l_3]$ and $W_k^{(M,N)}[l_1,l_2,0]$. The range of
$(a,c)$ for $W_k^{(M,N)}[l_1,l_2,l_3]$ is
\bean
0\leq a\leq l_3,0\leq c\leq l_2-a.
\eean
The subset corresponding to $a=0$ is exactly the range for
$W_k^{(M,N)}[l_1,l_2,0]$.

Similarly, we have
\bea
l_1'=l_1-1+\bar c-\bar a,l_2'=k-\bar c,
\label{RANGE2}
\ena
for $W_k^{(M,N)}[l_1-1,l_2-1,l_3-1]$ where
\bean
0\leq\bar a\leq l_3-1,0\leq c\leq l_2-1-\bar a.
\eean
By the correspondence
\bean
\bar a=a-1,\bar c=c,
\eean
the range of (\ref{RANGE2}) overlaps the range of (\ref{RANGE}) exactly when
\bean
1\leq a\leq l_3,0\leq c\leq l_2-a.
\eean
\end{proof}

\begin{corollary}\label{TRIPLES0}
The triples
\bea
&&0\rightarrow \bV_k[l_1-1,l_2-1]\rightarrow \bV_k[l_1,l_2]
\rightarrow V_k[l_1,l_2]\rightarrow0\quad(l_1+l_2\leq k),\label{triple1}\\
&&0\rightarrow \bU_k[l_1-1,l_2-1]\rightarrow \bU_k[l_1,l_2]
\rightarrow U_k[l_1,l_2]\rightarrow0\quad(l_1+l_2\leq k),\label{triple4}\\
&&0\rightarrow \bW_k[k-l_2-1,k-l_1-1]\rightarrow
\bW_k[l_1,l_2]\rightarrow W_k[l_1,l_2]\rightarrow0,\quad(l_1+l_2\geq k)\label{triple3}
\ena
are exact.
\end{corollary}
\begin{proof}
The exactness of triple \Ref{triple3} follows from Propositions
\ref{triple prop2} and \ref{tr prop}. The other two cases are obtained from
\Ref{triple6} by applying the automorphisms $T^{-1}$ and $\iota T^{-1}$.
\end{proof}

\subsection{Injectivity of coproduct}\label{multiplication}
Recall that the coproduct $\Delta$ in the universal enveloping algebra of 
Lie algebra $\widehat{\mathfrak{sl}}_2$ is defined by the condition
$\Delta(g)=g\otimes1+1\otimes g$ for $g\in\widehat{\mathfrak{sl}}_2$.

There exists an $\widehat{\mathfrak{sl}}_2$-linear map
\bean
\Delta:L_{k^{(1)}+k^{(2)},l^{(1)}+l^{(2)}}
&\rightarrow &L_{k^{(1)},l^{(1)}}\otimes L_{k^{(2)},l^{(2)}},\\
\Delta(v_{k^{(1)}+k^{(2)}}[l^{(1)}+l^{(2)}])&=&
v_{k^{(1)}}[l^{(1)}]\otimes v_{k^{(2)}}[l^{(2)}].
\eean
This map is obviously injective because of the irreducibility of 
$L_{k^{(1)}+k^{(2)},l^{(1)}+l^{(2)}}$. Similar mappings also exist for the
$\cH_k$-modules. However, the injectivity is not clear (and non-trivial)
because they are not
irreducible modules. In this section, we prove the injectivity of such maps
for $W_k[l_1,l_2,l_3]$. In fact we will show the injectivity of the coproduct
for the spaces $W_k^{(M,N)}[l_1,l_2,l_3]$. As a corollary, we obtain the
injectivity of the map
\bean
L_{k^{(1)}+k^{(2)},l^{(1)}+l^{(2)}}^{(M,N)}&\rightarrow
L_{k^{(1)},l^{(1)}}^{(M,N)}\otimes L_{k^{(2)},l^{(2)}}^{(M,N)}.
\eean

\begin{remark}
Recall the vector $u_k[l]=\frac{f_0^l}{l!}v_k[l]$ given by (\ref{UVEC}).
The above map $\Delta$
is equivalently characterized by the condition
$$
\Delta(u_{k^{(1)}+k^{(2)}}[l^{(1)}+l^{(2)}])=
u_{k^{(1)}}[l^{(1)}]\otimes u_{k^{(2)}}[l^{(2)}].
$$
\end{remark}

Define the coproduct $\Delta$ in the universal enveloping algebra of 
the affine Heisenberg algebra $\cH$ by the condition
$\Delta(g)=g\otimes1+1\otimes g$ for $g\in\cH$.

Let $k^{(j)},l_i^{(j)}$ $(i=1,2,3,\ j=1,2)$ be such that $k^{(1)}+k^{(2)}=k$
and $l_i^{(1)}+l_i^{(2)}=l_i$ $(1\le i\le3)$. Let
\bean
\Delta&:&W_k[l_1,l_2,l_3]\rightarrow
W_{k^{(1)}}[l^{(1)}_1,l^{(1)}_2,l^{(1)}_3]\otimes
W_{k^{(2)}}[l^{(2)}_1,l^{(2)}_2,l^{(2)}_3]
\eean
be the $\cH$-linear map such that
$\Delta(v)=v^{(1)}\otimes v^{(2)}$, where $v$, $v^{(1)}$ and $v^{(2)}$
are highest weight vectors of $W_k[l_1,l_2,l_3]$,
$W_{k^{(1)}}[l^{(1)}_1,l^{(1)}_2,l^{(1)}_3]$ and
$W_{k^{(2)}}[l^{(2)}_1,l^{(2)}_2,l^{(2)}_3]$ respectively.

The coproduct $\Delta$ descends to the map $\Delta^{(M,N)}$ on the coinvariants
\bean
\Delta^{(M,N)}:W_k^{(M,N)}[l_1,l_2,l_3]\rightarrow
W_{k^{(1)}}^{(M,N)}[l^{(1)}_1,l^{(1)}_2,l^{(1)}_3]\otimes
W_{k^{(2)}}^{(M,N)}[l^{(2)}_1,l^{(2)}_2,l^{(2)}_3].
\eean
Note that the cutoff (\ref{MD2}) does not break the well-definedness
of the map.

For simplicity, we suppress the dependence of these mappings on
$k^{(i)}$, $l^{(j)}_i$ in our notation,.

We start from
\begin{lemma}\label{SIMPLE}
Suppose that we have filtrations of vectors spaces:
$$
W^{(i)}=W^{(i)}_0\supset W^{(i)}_1\supset
\cdots\supset W^{(i)}_m\supset W^{(i)}_{m+1}=0\quad(i=1,2).
$$
We have the filtration of $W=W^{(1)}\otimes W^{(2)}$ consisting of the
subspaces $W_n=\sum_{a+b=n}W^{(1)}_a\otimes W^{(2)}_b$.
Then, the mapping
\begin{equation}\label{WELL}
W_n\rightarrow\oplus_{a+b=n}{\rm Gr}_a(W^{(1)})\otimes{\rm Gr}_b(W^{(2)})
\end{equation}
induced from the canonical surjections
$$
\pi_{a,b}:W^{(1)}_a\otimes W^{(2)}_b\rightarrow
{\rm Gr}_a(W^{(1)})\otimes{\rm Gr}_b(W^{(2)})
$$
is well-defined and induces the isomorphism
$$
{\rm Gr}_n(W)\rightarrow
\oplus_{a+b=n}{\rm Gr}_a(W^{(1)})\otimes{\rm Gr}_b(W^{(2)}).
$$
\end{lemma}
\begin{proof}
Suppose that an element $w\in W_n$ can be expressed in two different ways:
$w=\sum_ix_i=\sum_iy_i$ where $x_i,y_i\in W^{(1)}_i\otimes W^{(2)}_{n-i}$.
We must show that $\pi_{i,n-i}(x_i-y_i)=0$ for each $i$.

Consider a pair of vectors spaces and their subspaces:
$V\supset V_1,U\supset U_1$. Then, we have
$(V\otimes U_1)\cap(V_1\otimes U)=V_1\otimes U_1$.

Using this one can prove that
$$
x_i-y_i\in W^{(1)}_i\otimes W^{(2)}_{n-i+1}+
 W^{(1)}_{i+1}\otimes W^{(2)}_{n-i}.
$$
The welldefinedness of (\ref{WELL}) follows from this, and the rest of
the lemma is easy.
\end{proof}

\begin{proposition}\label{INJECTIVE}
The mapping $\Delta^{(M,N)}$ is injective.
\end{proposition}
\begin{proof}
It is enough to prove the case $k^{(1)}=k-1$ and $k^{(2)}=1$.
There are four cases:

$(i)$ $\Delta^{(M,N)}:W_k^{(M,N)}[l_1,l_2,l_3]
\rightarrow W_{k-1}^{(M,N)}[l_1-1,l_2,l_3]\otimes W_1^{(M,N)}[1,0,0]$,

$(ii)$ $\Delta^{(M,N)}:W_k^{(M,N)}[l_1,l_2,l_3]
\rightarrow W_{k-1}^{(M,N)}[l_1,l_2-1,l_3]\otimes W_1^{(M,N)}[0,1,0]$,

$(iii)$ $\Delta^{(M,N)}:W_k^{(M,N)}[l_1,l_2,l_3]
\rightarrow W_{k-1}^{(M,N)}[l_1-1,l_2-1,l_3]\otimes W_1^{(M,N)}[1,1,0]$,

$(iv)$ $\Delta^{(M,N)}:W_k^{(M,N)}[l_1,l_2,l_3]
\rightarrow W_{k-1}^{(M,N)}[l_1-1,l_2-1,l_3-1]\otimes W_1^{(M,N)}[1,1,1]$.

We use induction on $(M,N)$. If $(M,N)=(0,0)$,
the assertion is clear from (\ref{INIT}).

In the induction, we use the recurrence relations for the space
$W_1^{(M,N)}[l_1,l_2,l_3]$ given by Theorem \ref{R}.
We have the following bijections:
\bean
W_1^{(M,N-1)}[1,1]&\rightarrow&{\rm Gr}_{0,0}(W_1^{(M,N)}[1,0,0]),\\
W_1^{(M,N-1)}[0,1]&\rightarrow&{\rm Gr}_{0,0}(W_1^{(M,N)}[0,1,0]),\\
W_1^{(M,N-1)}[1,0]&\rightarrow&{\rm Gr}_{0,1}(W_1^{(M,N)}[0,1,0]),\\
W_1^{(M,N-1)}[1,1]&\rightarrow&{\rm Gr}_{0,0}(W_1^{(M,N)}[1,1,0]),\\
W_1^{(M,N-1)}[1,0]&\rightarrow&{\rm Gr}_{0,1}(W_1^{(M,N)}[1,1,0]),\\
W_1^{(M,N-1)}[1,1]&\rightarrow&{\rm Gr}_{0,0}(W_1^{(M,N)}[1,1,1]),\\
W_1^{(M,N-1)}[1,0]&\rightarrow&{\rm Gr}_{0,1}(W_1^{(M,N)}[1,1,1]),\\
W_1^{(M,N-1)}[0,1]&\rightarrow&{\rm Gr}_{1,0}(W_1^{(M,N)}[1,1,1]).
\eean
Let us prove the injectivity of $(i)$. We use the following abbreviated
notations. $W=W_k^{(M,N)}[l_1,l_2,l_3]$,
$W^{(1)}=W_{k-1}^{(M,N)}[l_1-1,l_2,l_3]$,
$W^{(2)}=W_1^{(M,N)}[1,0,0]$. We denote the highest weight vectors
in $W,W^{(1)},W^{(2)}$, respectively, by $v,v^{(1)},v^{(2)}$.
Note that $h_1v^{(2)}=f_1v^{(2)}=0$, and hence
\begin{equation}\label{ONETERM}
\Delta(h_1^af_1^c)(v^{(1)}\otimes v^{(2)})=(h_1^af_1^cv^{(1)})\otimes v^{(2)}.
\end{equation}
Consider the following mapping diagram.
\[
\xymatrix{
C_{a,c}(W)\ar[r]\ar[d]
&C_{a,c}(W^{(1)})\otimes C_{0,0}(W^{(2)})\ar[d]
\\
{\rm Gr}_{a,c}(W)\ar[r]
&{\rm Gr}_{a,c}(W^{(1)})\otimes{\rm Gr}_{0,0}(W^{(2)})
\\
W^{(M,N-1)}[l'_1,l'_2]\ar[r]\ar[u]
&W^{(M,N-1)}[l'_1-1,l'_2-1]\otimes W_1^{(M,N-1)}[1,1]\ar[u]
}
\]
The first and the third horizontal arrows are the coproduct maps.
The downward vertical arrows are canonical surjections, and the
upward vertical arrows are the canonical bijections of Proposition
\ref{RECURSION}.
The second horizontal arrow is induced from either of the other horizontal
arrows. Noting that the shift automorphism $T$ and the coproduct $\Delta$
are commutative, one can check that the induced map is unique, and thereby
we have a commutative diagram.

By the induction hypothesis, the third horizontal arrow
is injective. Therefore, the second horizontal arrow is also injective. From
this follows that the kernel of the first horizontal arrow is included
in a smaller subset in the filtration. Repeating this argument, i.e.,
using induction on $(a,c)$, one can show the injectivity
of the coproduct map $\Delta^{(M,N)}$.

In other three cases, the relation (\ref{ONETERM}) is modified.
Let us consider the case (ii).
Using similar abbreviated notations, we have
\begin{equation}\label{TWOTERM}
\Delta(h_1^af_1^c)(v^{(1)}\otimes v^{(2)})=(h_1^af_1^cv^{(1)})\otimes v^{(2)}
+c(h_1^af_1^{c-1}v^{(1)})\otimes f_1v^{(2)}.
\end{equation}
This induces a mapping of the form
\begin{equation}
C_{a,c}(W)\rightarrow C_{a,c}(W^{(1)})\otimes C_{0,0}(W^{(2)})
+C_{a,c-1}(W^{(1)})\otimes C_{0,1}(W^{(2)}).
\end{equation}
Here we use
$W=W_k^{(M,N)}[l_1,l_2,l_3],
W^{(1)}=W_{k-1}^{(M,N)}[l_1,l_2-1,l_3],
W^{(2)}=W_1^{(M,N)}[0,1,0]$. 

If $c\not=0$, we consider the following mapping diagram.
\[
\xymatrix{
C_{a,c}(W)\ar[r]\ar[d]
&C_{a,c}(W^{(1)})\otimes C_{0,0}(W^{(2)})
+C_{a,c-1}(W^{(1)})\otimes C_{0,1}(W^{(2)})\ar[d]
\\
{\rm Gr}_{a,c}(W)\ar[r]
&{\rm Gr}_{a,c-1}(W^{(1)})\otimes{\rm Gr}_{0,1}(W^{(2)})
\\
W^{(M,N-1)}[l'_1,l'_2]\ar[r]\ar[u]
&W^{(M,N-1)}[l'_1-1,l'_2]\otimes W_1^{(M,N-1)}[1,0]\ar[u]
}
\]
The first horizontal arrow is $\frac1c\Delta$. The downward vertical arrow
in the right is assured by Lemma \ref{SIMPLE}. The rest of proof is
similar.

If $c=0$ the second term in the right hand side of (\ref{TWOTERM}) is absent.
In this case, we use
\[
\xymatrix{
C_{a,c}(W)\ar[r]\ar[d]
&C_{a,c}(W^{(1)})\otimes C_{0,0}(W^{(2)})\ar[d]
\\
{\rm Gr}_{a,c}(W)\ar[r]
&{\rm Gr}_{a,c}(W^{(1)})\otimes{\rm Gr}_{0,0}(W^{(2)})
\\
W^{(M,N-1)}[l'_1,l'_2]\ar[r]\ar[u]
&W^{(M,N-1)}[l'_1,l'_2-1]\otimes W_1^{(M,N-1)}[0,1]\ar[u]
}
\]
The rest of the proof is similar.

In the cases (iii) and (iv), the proof will go with
${\rm Gr}_{a,c}(W^{(1)})\otimes{\rm Gr}_{0,0}(W^{(2)})$
using the bijection of the form
$$
W^{(M,N-1)}[l'_1-1,l'_2-1]\rightarrow{\rm Gr}_{a,c}(W^{(1)}).
$$
\end{proof}

Now we are ready to show the injectivity of coproduct for the
$\widehat{\mathfrak{sl}}_2$ spaces of coinvariants.

\begin{theorem}\label{WINJ}
The mapping
\bean
\Delta^{(M,N)}&:&L_{k^{(1)}+k^{(2)},l^{(1)}+l^{(2)}}^{(M,N)}\rightarrow
L_{k^{(1)},l^{(1)}}^{(M,N)}\otimes L_{k^{(2)},l^{(2)}}^{(M,N)}
\eean
is injective.
\end{theorem}
\begin{proof}
We assume that $N>0$. Consider the following commutative diagram.
\[
\xymatrix{
L^{(M,N)}_{k,l}\ar[r]&
L^{(M,N)}_{k^{(1)},l^{(1)}}\otimes L^{(M,N)}_{k^{(2)},l^{(2)}}\\
E_i(L_{k,l})\ar[u]^{\iota_{k,l}}
\ar[r]
\ar[d]^{\pi_{k,l}}
&
\sum_{i_1+i_2=i}E_{i_1}(L_{k^{(1)},l^{(1)}})\otimes
E_{i_2}(L_{k^{(2)},l^{(2)}})
\ar[u]^{\iota_{k^{(1)},l^{(1)}}\otimes\iota_{k^{(2)},l^{(2)}}}
\ar[d]^{\pi_{k^{(1)},l^{(1)}}\otimes\pi_{k^{(2)},l^{(2)}}}\\
{\rm Gr}^E_i(L_{k,l})^{(M,N)}\ar[r]&
\oplus_{i_1+i_2=i}{\rm Gr}^E_{i_1}(L_{k^{(1)},l^{(1)}})^{(M,N)}
\otimes{\rm Gr}^E_{i_2}(L_{k^{(2)},l^{(2)}})^{(M,N)}
}\]
The horizontal arrows are coproducts, the vertical up-arrows are
compositions of canonical injections and surjections,
and the vertical down-arrows are canonical
surjections. The only non-trivial one is
$$
\sum_{i_1+i_2=i}E_{i_1}(L_{k^{(1)},l^{(1)}})\otimes
E_{i_2}(L_{k^{(2)},l^{(2)}})
\rightarrow
\oplus_{i_1+i_2=i}{\rm Gr}^E_{i_1}(L_{k^{(1)},l^{(1)}})^{(M,N)}
\otimes{\rm Gr}^E_{i_2}(L_{k^{(2)},l^{(2)}})^{(M,N)}
$$
which follows from Lemma \ref{SIMPLE}.

We are to prove that the first horizontal arrow is injective.
We know that the third horizontal arrow is injective (Proposition
\ref{INJECTIVE} up to the automorphisms $T'$ given by (\ref{T'})).
Note also that because of (\ref{NOJUMP}) we can apply Remark \ref{rem}.

Suppose $w\in L^{(M,N)}_{k,l}$ belongs to the kernel of $\Delta^{(M,N)}$.
If $w\not=0$, we take the smallest values $i$ such that $w$ belongs to the
image of $E_i(L_{k,l})$. Take a preimage $w'\in E_i(L_{k,l})$ of $w$. We have
$i>0$ because $\Delta^{(M,N)}(u_k[l])\not=0$. We show a contradiction by
finding $w''\in E_{i-1}(L_{k,l})$ satisfying $w=\iota_{k,l}(w'')$.

For simplicity set
$\iota=\iota_{k^{(1)},l^{(1)}}\otimes\iota_{k^{(2)},l^{(2)}}$
and
$\pi=\pi_{k^{(1)},l^{(1)}}\otimes\pi_{k^{(2)},l^{(2)}}$.
Since $\Delta^{(M,N)}(w)=0$, $\Delta(w')$ belongs to ${\rm Ker}\,\iota$. 
We will show that 
\begin{equation}
{\rm Ker}\,\iota\subset{\rm Ker}\,\pi.\label{INCL}
\end{equation}

Prepare a set of vectors
$B^{(a)}_j=\{v^{(a)}_{j,\alpha}\in E_j(L_{k^{(a)},l^{(a)}})
;1\leq\alpha\leq{\rm dim}\,{\rm Gr}^E_j(L_{k^{(a)},l^{(a)}})^{(M,N)}\}$
for each $(a,j)$ such that
$\pi_{k^{(a)},l^{(a)}}(B^{(a)}_j)$ forms a basis of
${\rm Gr}^E_j(L_{k^{(a)},l^{(a)}})^{(M,N)}$. Because of Remark \ref{rem}
$\iota_{k^{(a)},l^{(a)}}(B^{(a)}_j)$ forms a basis of
$L^{(M,N)}_{k^{(a)},l^{(a)}}$.

Suppose that $v=\sum_{i_1+i_2=i}v_{i_1,i_2}$ where $v_{i_1,i_2}\in
E_{i_1}(L_{k^{(1)},l^{(1)}})\otimes E_{i_2}(L_{k^{(2)},l^{(2)}})$.
We write
$$
\pi(v_{i_1,i_2})=\sum_{\alpha,\beta}c^{\alpha,\beta}_{i_1,i_2}
\pi_{k^{(1)},l^{(1)}}(v^{(1)}_{i_1,\alpha})\otimes
\pi_{k^{(2)},l^{(2)}}(v^{(2)}_{i_2,\beta}).
$$
We have
$$
v_{i_1,i_2}-\sum_{\alpha,\beta}c^{\alpha,\beta}_{i_1,i_2}
v^{(1)}_{i_1,\alpha}\otimes v^{(2)}_{i_2,\beta}\in{\rm Ker}\,\pi
\cap
E_{i_1}(L_{k^{(1)},l^{(1)}})\otimes E_{i_2}(L_{k^{(2)},l^{(2)}}).
$$

Note that
\begin{eqnarray*}
&&{\rm Gr}^E_{i_a}(L_{k^{(a)},l^{(a)}})^{(M,N)}\\
&&\quad\simeq
E_{i_a}(L_{k^{(a},l^{(a)}})/\Bigl(
\sum_{j\geq M}e_jE_{i_a-1}(L_{k^{(a)},l^{(a)}})
+\sum_{j\geq N}f_jE_{i_a-1}(L_{k^{(a)},l^{(a)}})
+E_{i_a-1}(L_{k^{(a)},l^{(a)}})\Bigr).
\end{eqnarray*}
{}From this we conclude that
$v-\sum_{i_1+i_2=i\atop\alpha,\beta}
c^{\alpha,\beta}_{i_1,i_2}v^{(1)}_{i_1,\alpha}\otimes v^{(2)}_{i_2,\beta}$
belongs to
\begin{eqnarray*}
&&\sum_{i_1+i_2=i-1}
E_{i_1}(L_{k^{(1)},l^{(1)}})\otimes E_{i_2}(L_{k^{(2)},l^{(2)}})\\
&&
+\sum_{i_1+i_2=i}\Biggl(\sum_{j\geq M}
\Bigl(e_jE_{i_1-1}(L_{k^{(1)},l^{(1)}})\Bigr)
\otimes E_{i_2}(L_{k^{(2)},l^{(2)}})
+\sum_{j\geq N}\Bigl(f_jE_{i_1-1}(L_{k^{(1)},l^{(1)}})\Bigr)
\otimes E_{i_2}(L_{k^{(2)},l^{(2)}})\\
&&+\sum_{j\geq M}
E_{i_1}(L_{k^{(1)},l^{(1)}})
\otimes\Bigl(e_jE_{i_2-1}(L_{k^{(2)},l^{(2)}})\Bigr)
+\sum_{j\geq N}E_{i_1}(L_{k^{(1)},l^{(1)}})
\otimes\Bigl(f_jE_{i_2-1}(L_{k^{(2)},l^{(2)}})\Bigr)\Biggr).
\end{eqnarray*}
Repeating the argument for $E_{i-1}(L_{k,l})$ and so on, we obtain
$$
v=\sum_{i_1+i_2\leq i\atop \alpha,\beta}
c^{\alpha,\beta}_{i_1,i_2}
v^{(1)}_{i_1,\alpha}\otimes v^{(2)}_{i_2,\beta}
\bmod{\rm Ker}\,\iota\cap{\rm Ker}\,\pi.
$$

If $v\in{\rm Ker}\,\iota$, thus we have 
$c^{\alpha,\beta}_{i_1,i_2}=0$ for all $\alpha,\beta,i_1,i_2$,
and thereby $v\in{\rm Ker}\,\pi$.

We return to $w'\in E_i(L_{k,l})$ which satisfies
$\Delta(w')\in{\rm Ker}\,\pi$. Using (\ref{INCL}) and the injectivity of
the third horizontal arrow, we conclude that $\pi_{k,l}(w')=0$.

This implies
$$
w'\in\sum_{j\geq M}e_jE_{i-1}(L_{k,l})
+\sum_{j\geq N}f_jE_{i-1}(L_{k,l})
+E_{i-1}(L_{k,l}).
$$
Therefore, we can replace $w'$ by $w''\in E_{i-1}(L_{k,l})$
keeping the property $\iota(w'')=w$.
\end{proof}
\subsection{Character identities}\label{char identities section}
In this section we collect the relations between characters of coinvariants of different modules. All these equations are easy corollaries of the results proved in the paper.

Recall that for a $\cH_k$ module $V$ we have a
triple grading such that
\bean
{\rm deg}\, v=(0,0,0),\quad {\rm deg}\,e_i=(1,0,i),\quad 
{\rm deg}\,f_i=(0,1,i),\quad
{\rm deg}\,h_i=(1,1,i),
\eean
where $v$ is the highest weight vector in $V$.
Let as before $V_{m,n,e}$ be the subspace of degree $(m,n,e)$ and 
$V_{m,n}=\oplus_eV_{m,n,e}$.

Define the character of $V$ by 
\be
\chi(V)(q,z_1,z_2)=\sum_{m,n,e\in\Z_{\geq 0}}z_1^mz_2^nq^e\dim V_{m,n,e}.
\ee

For $M,N\geq 0$, the space of coinvariants $V^{(M,N)}$ is a quotient of $V$ by a graded ideal. 
Therefore, we have a well-defined induced character
of the coinvariant space $\chi(V^{(M,N)})(q,z_1,z_2)$. 

In the limit $M,N\to \infty$ we obviously obtain the character of the whole
module,
\be
\chi(V^{(\infty,\infty)})(q,z_1,z_2)=\chi(V)(q,z_1,z_2).
\ee

Below we assume that $M,N\geq 0$, moreover, $N>0$ for the relations which
include $V_k^{(M,N)}[l_1,l_2,l_3]$  and $M>0$ for the relations which include
$U_k^{(M,N)}[l_1,l_2,l_3]$

First, we give the relation between the characters of the Heisenberg and
$\wsl$ coinvariants. 

We have (see Corollaries \ref{COINV COL}, \ref{GRA})
\bea
z^l\chi(V_k[k-l,l]^{(M,N)})(q,z^2,z^{-2})
=\chi_{k,l}^{(M,N)} \qquad (M>0,N\geq 0),\label{V=L}\\
z^{-l} \chi(U_k[l,k-l]^{(M,N)})(q,z^2,z^{-2})
=\chi_{k,l}^{(M,N)} \qquad (M\geq 0,N>0).\label{U=L}
\ena
In particular,
\be
\chi(V_k[k-l,l]^{(M,N)}(1,1,1)=d^{(M+N)}_{k,l},\qquad
 \chi(U_k[l,k-l]^{(M,N)}(1,1,1)=d^{(M+N)}_{k,l}.
\ee

By Lemma \ref{SPLIT} and Corollary \ref{SES cor}, the spaces $W$ refine these
identities, 
\bea\label{W=L}
\sum_{i=0}^lz_2^i\chi(W^{(M,N)}_k[k-l+i,k-i])(q,z_1,z_2)
=\chi(V_k[k-l,l]^{(M,N)})(q,z_1,z_2).
\ena
Also, we have
\be
z_2^i\chi(W^{(0,N)}_k[k-l+i,k-i])(q,z_1,z_2)=\chi_{k,l}^N[i,*](q,z_1,z_2)
\ee
(see \Ref{partial char} for the definition of the right hand side).

{}From \Ref{V=W}, \Ref{V=U}, we have 
\bean
\chi(W_k[l_1,l_2,l_3]^{(M,N)})(q,z_1,z_2)&=&\chi(V_k[l_1,l_2,l_3]^{(M,N-1)})(q,z_1,qz_2),\\
\chi(U_k[l_1,l_2,l_3]^{(M,N)})(q,z_1,z_2)&=&\chi(V_k[l_2,l_1,l_3]^{(N,M)})(q,z_2,z_1).
\eean

The short exact sequences \Ref{triple2}, \Ref{triple5}, \Ref{triple1},
\Ref{triple4}, \Ref{triple3} and \Ref{triple6} give us the following relations.
\bean 
\chi(V_k[l_1,l_2]^{(M,N)})&=&\chi(\bV_k[l_1,l_2]^{(M,N)}
)-z_1z_2\chi(\bV_k[l_1-1,l_2-1]^{(M,N)}),\\
\chi(U_k[l_1,l_2]^{(M,N)})&=&\chi(\bU_k[l_1,l_2]^{(M,N)})-z_1z_2\chi(\bU_k[l_1-1,l_2-1]^{(M,N)}),\\
\chi(W_k[l_1,l_2]^{(M,N)})&=&\chi(\bW_k[l_1,l_2]^{(M,N)})-(qz_1z_2)^{l_1+l_2-k+1}\chi(\bW_k[k-l_2-1,k-l_1-1]^{(M,N)})
\eean
and
\bean
\chi(V_k[l_1,l_2]^{(M,N)})&=&\chi(V_k[l_1,l_2-1]^{(M,N)})+z_2^{l_2}\chi(W_k[l_1+l_2,k-l_2]^{(M,N)}),\\
\chi(U_k[l_1,l_2]^{(M,N)})&=&\chi(U_k[l_1-1,l_2-1]^{(M,N)})+z_1^{l_1}\chi(W_k[k-l_1,l_1+l_2]^{(M,N)}),\\
\chi(W_k[l_1,l_2,l_3]^{(M,N)})&=&\chi(W_k[l_1,l_2,0]^{(M,N)})+qz_1z_2\chi(W_k[l_1-1,l_2-1,l_3-1]^{(M,N)}).
\eean

The main recursion (see Theorem \ref{R}) gives us the following relation
\be
\chi(W_k[l_1,l_2,l_3]^{(M,N)})(q,z_1,z_2)=\sum_{\stackrel{0\leq a\leq
l_3}{0\leq c\leq l_2-a}} q^{a+c}z_1^a
z_2^{a+c}\chi(W_k[l_1',l_2']^{(M,N-1)})(q,z_1,qz_2). 
\ee
In the case $M=0$, $l_3=l_1+l_2-k$, this relation coincides with \Ref{shift recursion}.

\section{Appendix}
In this appendix we give several technical lemmas about
some vectors being zero in the $\cH_k$-modules.
\subsection{Zero vector lemmas}
Let $V$ be a $\cH_k$-module and $u\in V$ be a vector.

\begin{lemma}\label{LM0}
If $e_1^{a+1}u=f_{-1}u=0$, then 
\bean
e_1^{a-c}h_0^{c+1}u=0,\qquad 0\leq c\leq a.
\eean 
If $f_0^{b+1}u=e_0u=0$, then
\bean
f_0^{b-c}h_0^{c+1}u=0,\qquad 0\leq c\leq b.
\eean
\end{lemma}
\begin{proof}
Since $({\rm ad}f_{-1})^{c}(e_1^{a+1})u=0$,  by
using $[e_1,f_{-1}]=h_0$ and $[h_0,f_{-1}]=0$, we obtain
the first equation of the lemma. The proof of second one is similar.
\end{proof}

\begin{lemma}\label{LM1}
If $e_1^au=0$, then
\bean
e_1^{a+b}f_1^bu&=&0.
\label{L1}
\eean
\end{lemma}
\begin{proof}
First, we show 
\bea
e_1^{a+b-1}f_1^be_1u&=&0
\label{INDHYP}
\ena
by induction on $b$. The case $b=0$ is the assumption of the Lemma.
Suppose (\ref{INDHYP}) is true for $b-1$. Then
\bean
e_1^{a+b-1}f_1^be_1u&=&e_1[e_1^{a+b-2},f_1]f_1^{b-1}e_1u\nonumber\\
&=&(a+b-2)h_2e_1^{a+b-2}f_1^{b-1}e_1u\nonumber\\
&=&0,
\eean
and \Ref{INDHYP} holds for $b$.
Next, we have
\bean
e_1^{a+b}f_1^bu&=&e_1^{a+b-1}[e_1,f_1^b]u\\
&=&bh_2e_1^{a+b-1}f_1^{b-1}u.
\eean
Therefore, the lemma follows by induction on $b$.
\end{proof}

\begin{lemma}\label{LM2}
If $e_1^au=f_0u=0$, then
\bean
e_1^{a-b+c}h_1^bf_1^cu&=&0.
\eean
\end{lemma}
\begin{proof}
We use induction on $b$. The case $b=0$ follows from Lemma \ref{LM1}.
We have
\bean
e_1^{a-b+c}h_1^bf_1^cu&=&e_1^{a-b+c}h_1^{b-1}[e_1,f_0]f_1^cu\\
&=&-e_1^{a-b+c}h_1^{b-1}f_0e_1f_1^cu.
\eean
Using the induction hypothesis, we have
\bean
e_1^{a-b+c}h_1^bf_1^cu&=&[f_0,e_1^{a-b+c}]h_1^{b-1}e_1f_1^cu\\
&=&-(a-b+c)e_1^{a-b+c}h_1^bf_1^cu.
\eean
The assertion follows from this.
\end{proof}

\begin{lemma}\label{LM3}
If $e_1^au=h_1^bf_1^{c+1}u=f_0u=0$, then
\bean
h_2^{a-b}h_1^bf_1^cu&=&0.
\eean
\end{lemma}
\begin{proof}
First, we prove
\bean
h_1^bf_1^{c+1+n}e_1^nu&=&0
\eean
by induction on $n$. In fact, 
\bean
h_1^bf_1^{c+1+n}e_1^nu&=&[h_1^bf_1^{c+1+n},e_1]e_1^{n-1}u\\
&=&-(c+1+n)h_2h_1^bf_1^{c+n}e_1^{n-1}u\\
&=&0.
\eean
Next, we proceed as follows:
\bean
h_2h_1^bf_1^{c+n}e_1^nu&=&h_1^bf_1^{c+n}[e_1,f_1]e_1^nu\\
&=&[h_1^bf_1^{c+n},e_1]f_1e_1^nu-h_1^bf_1^{c+n+1}e_1^{n+1}u\\
&=&-(c+n)h_2h_1^bf_1^{c+n}e_1^nu
-h_1^bf_1^{c+n+1}e_1^{n+1}u.
\eean
By repeating this argument, the assertion reduces to showing that
\bean
h_1^bf_1^{c+a-b}e_1^{a-b}u=0.
\eean
This follows from Lemma \ref{LM2} with $c=0$.
\end{proof}


\begin{thebibliography}{FKLMM}
\bibitem[B]{B} R.J. Baxter, {\it Exactly solvable models in
Statistical mechanics}, Academic Press, London, 1982.

\bibitem[DJKMO]{DJKMO} E. Date, M. Jimbo, A. Kuniba, T. Miwa, and
M. Okado, {\it One dimensional configuration sums in vertex models and
affine Lie algebra characters}, Lett. Math. Phys., 17 (1989), 69-77.

\bibitem[FL]{FL}B. Feigin and S. Loktev, {\it On generalized Kostka polynomials
and quantum Verlinde rule}, math.QA/9812093.

\bibitem[FM]{FM}B. Feigin and T. Miwa,{\it Extended vertex operator algebras
and monomial bases}, math.QA/9901067.  

\bibitem[FS]{FS}B. Feigin and A. Stoyanovsky, {\it Quasi-particles
models for the representations of Lie algebras and geometry of flag
manifold}, hep-th/9308079, RIMS 942; {\it Functional models for the
representations of current algebras and the semi-infinite Schubert
cells}, Funct. Anal. Appl. {\bf 28} (1994), 55--72.

\bibitem[Fi]{Fi}M. Finkelberg, {\it An equivalence of fusion categories},
Geom. Funct. Anal. {\bf 6}(1996)249-267.

\bibitem[FKLMM]{FKLMM} B. Feigin, R. Kedem, S. Loktev, T. Miwa,
E. Mukhin, {\it Combinatorics of the $\wsl$ spaces of coinvariants},
math-ph/9908003, RIMS 1243.

\bibitem[KKMM]{KKMM} R. Kedem, T. Klassen, B. McCoy, E. Melzer,
{\it Fermionic sum representations for conformal field theory
characters}. Phys. Lett.  {\bf B 307} (1993), 68--76.

\bibitem[S]{S}G. Segal, Geometric aspect of quantum field theory,
in Proceedings of the ICM at Kyoto, 1990.

\bibitem[St]{St}A. Stoyanovsky, {\it Lie algebra deformation and character formulas},
Func. Anal. Its Appl.,{\bf 32}(1998)66-68.

\bibitem[TUY]{TUY}A. Tsuchiya, K. Ueno and Y. Yamada,{\it Conformal field
theory on the universal family of stable curves with gauge symmetry},
Adv. Stud. Pure Math., {\bf 19}(1989)459-466.
\end{thebibliography}
\end{document}